\documentclass[12pt,a4paper,reqno]{amsart}
\usepackage{amssymb}
\usepackage{dcpic,pictex}
\usepackage{multirow}
\usepackage{graphicx}
\usepackage{slashbox}
\usepackage{cite}

\usepackage{hyperref}
\usepackage[main=english,french]{babel}
\newtheorem{definition}{Definition}[section]

\newtheorem{theorem}{Theorem}[section]
\newtheorem{remark}{Remark}[section]
\newcommand{\D}{\mathbb{D}}

\newcommand{\I}{\mathbb{I}}
\newcommand{\R}{\mathbb{R}}
\newcommand{\C}{\mathbb{C}}
\newcommand{\N}{\mathbb{N}}

\begin{document}


\title[Operational Calculus]{Operational Calculus for the general 
       fractional  derivatives with the Sonine kernels}

\author{Yuri Luchko}
\curraddr{Beuth Technical University of Applied Sciences Berlin,  
     Department of  Mathematics, Physics, and Chemistry,  
     Luxemburger Str. 10,  
     13353 Berlin,
Germany}
\email{luchko@beuth-hochschule.de}

\dedicatory{}

\begin{abstract}
In this paper, we first address the general fractional integrals and derivatives with the Sonine kernels that possess the integrable singularities of power function type at the point zero. Both particular cases and compositions of these operators are discussed. Then we proceed with a construction of an operational calculus of the Mikusi\'nski type for the general fractional derivatives with the Sonine kernels. This operational calculus is applied for analytical treatment of some initial value problems for the fractional differential equations with the general fractional derivatives. The solutions are expressed in form of the convolution series that generalize the power  series for the exponential and the Mittag-Leffler functions. 
\end{abstract}

\maketitle

\section{Introduction}
\label{sec:1}

\setcounter{section}{1}
\setcounter{equation}{0}\setcounter{theorem}{0}

The evolution equations containing the integro-differential operators of convolution type
\begin{equation}
\label{FDR-L} 
(\D_{(k)}\, f)(t) = \frac{d}{dt}\, \int_0^t k(t-\tau)f(\tau)\, d\tau,\ t>0
\end{equation}
were studied in the framework of the abstract Volterra integral equations on the Banach spaces for more than a half century. In particular, in \cite{Cle84}, the abstract Volterra integral equations with the operator \eqref{FDR-L} were addressed in the case of the completely  positive kernels $k\in L^1(0,+\infty)$. A completely  positive kernel is defined as the one that satisfies the relation
\begin{equation}
\label{Cle_1} 
a\, k(t) \, +\, \int_0^t k(t-\tau)l(\tau)\, d\tau\, = 1,\  t>0
\end{equation}
with $a \ge 0$ and $l\in L^1(0,+\infty)$, the function $l$ being non-negative and non-increasing. 

However, until recently, no interpretation of these general results in the framework of the Fractional Calculus (FC) was suggested. The starting point for involving the general integro-differential operators of convolution type into FC research was the  paper \cite{Koch11} (see also \cite{Han20,Luc21,LucYam20}). In \cite{Koch11}, the operator \eqref{FDR-L} with the suitably defined kernels was called the general fractional derivative (GFD) of the Riemann-Liouville type.  The GFD of the Caputo type was defined in the form
\begin{equation}
\label{FDC}
( _*\D_{(k)}\, f) (t) =  (\D_{(k)}\, f) (t) - f(0)k(t).
\end{equation}

In \cite{Koch11}, a class of the kernels which allow an interpretation of the operators  \eqref{FDC} as a kind of the fractional derivatives was introduced. These are the kernels that satisfy the conditions K1)-K4) below:

\vspace{0.1cm}

\noindent
K1) The Laplace transform $\tilde k$ of $k$,
\begin{equation}
\label{Laplace} 
\tilde k(p) = ({\mathcal L}\, k)(p)\ =\ \int_0^{+\infty} k(t)\, e^{-pt}\, dt
\end{equation}
exists for all $p>0$,

\vspace{0.1cm}

\noindent
K2) $\tilde k(p)$ is a Stieltjes function (see \cite{[SSV]} for details regarding the Stieltjes functions),

\vspace{0.1cm}

\noindent
K3) $\tilde k(p) \to 0$ and $p \tilde k(p) \to +\infty$ as $p \to +\infty$,

\vspace{0.1cm}

\noindent
K4)  $\tilde k(p) \to +\infty$ and $p \tilde k(p) \to 0$ as $p \to 0$.

\vspace{0.1cm} 
In what follows, we denote the class of the Sonine kernels that satisfy the conditions K1)-K4) by $\mathcal{K}$. The operators \eqref{FDC} with the kernels from  $\mathcal{K}$ are left inverse to the appropriately defined fractional integrals, i.e., they satisfy the 1st Fundamental Theorem of FC 
(\cite{Koch11,Luc20,Luc21}). Moreover, the solutions to the time-fractional differential equations with the GFDs possess some typical features of the solutions to the evolution equations. More precisely, the statements (A) and (B) formulated below hold true:

\vspace{0.1cm}

\noindent
(A) For any $\lambda <0$, the initial value problem for the fractional relaxation equation 
\begin{equation} 
\label{relax}
\begin{cases}
( _*\D_{(k)} u) (t) = \lambda u(t),\   t >0, \\
u(0) =1
\end{cases}
\end{equation}
has a unique solution $u_\lambda = u_\lambda(t)$ that is continuous on $[0,\, +\infty)$ and infinitely differentiable and completely monotone on $\R_+$, i.e.,
\begin{equation} 
\label{cmf}
(-1)^n u_\lambda^{(n)}(t) \ge 0,\   t>0,\ n = 0,1,2\dots.
\end{equation}

\vspace{0.1cm}

\noindent
(B) The fundamental solution to the Cauchy problem for
the general time-fractional diffusion equation 
\begin{equation}
\label{eq_K}
\begin{cases}
( _*\D_{(k)}  u(x,\cdot))(t) \, =\, \Delta u(x,t),\ t>0,\ x\in \R^n,\\
u(x,0) = u_0(x),\ x\in \R^n
\end{cases}  
\end{equation}
is  locally integrable in $t$, infinitely differentiable for $x\not = 0$, and can be interpreted as a spatial probability density function evolving in time.     

Later on,  a maximum principle for the initial-boundary-value problems for the general time-fractional diffusion equation in form \eqref{eq_K} with a general second order spatial differential operator instead of the Laplace operator $\Delta$ was deduced in \cite{LucYam16}. This maximum principle was in particular  applied to prove solution uniqueness and for derivation of its a priory norm estimates. 

In the meantime, some advanced results for the ordinary and partial time-fractional  differential equations with the GFDs have been derived (see  \cite{Koch11,Koch19_2,KK,Sin18} for the ODEs and  \cite{Koch11,Koch19_1,Koch19_2,LucYam16,LucYam20,Sin20} for the PDEs with the GFDs). In \cite{JK17,KJ19_1,KJ19}, several inverse problems for the fractional differential equations with the GFDs were formulated and solved. 

In this paper, we address one of the powerful techniques for analytical treatment of the fractional differential equations that has been already successfully employed for equations with different kinds of the fractional derivatives. It is an operational calculus of the Mikusi\'nski type that we construct  for the GFD \eqref{FDC} with a Sonine kernel $k=k(t)$ that possesses an integrable singularity of power function type at the point zero. Then we provide some applications of this calculus to the fractional differential equations with the GFDs.  

Until now, operational calculi of the Mikusi\'nski type were constructed for the Riemann-Liouville fractional derivative (\cite{HadLuc}, \cite{Luc93}, \cite{LucSri95}),  the Caputo fractional derivative (\cite{LucGor99}), the  Hilfer fractional derivative  
(\cite{HLT09}), the  Caputo-type modification of the  
Erd\'elyi-Kober fractional derivative (\cite{Han}), and  the finite compositions  of the  Erd\'elyi-Kober fractional derivatives (\cite{BasLuc95, Luc93, LucYak94, YakLuc94}). 
In the publications mentioned above, these calculi were employed for analytical treatment of the initial-value problems for the fractional differential equations with different kinds of the fractional derivatives. We also refer to \cite{Luc99} and \cite{Luc19} for the surveys of the main results regarding the operational calculi for the fractional differential operators. In \cite{GorLuc97} and \cite{GLS97}, the Mikusi\'nski type operational calculi where employed for solving the generalized Abel integral equation and the integral equations with the Gauss hypergeometric function in the kernel, respectively. Thus, the operational calculi of he Mikusi\'nski type  are useful for analytical treatment of both differential and integral equations. 

The rest of the paper is organized as follows. In the second
section, the rings of functions that were employed so far for construction of the operational calculi of the Mikusi\'nski type are shortly discussed. 
The third section addresses the general fractional integrals (GFI) and derivatives with the Sonine kernels that possess the integrable singularities of power function type at the point zero.  The fourth section is devoted to construction of an operational calculus of the Mikusi\'nski type  for the GFD \eqref{FDC} with a Sonine kernel from this class. In the field of convolution quotients, both the GFIs and the GFDs are represented as multiplications with certain elements of this field. Moreover, we derive the basic operational relations that represent certain elements of the field of convolution quotients as conventional functions in form of the convolution series that generalize the power series for the exponential and the Mittag-Leffler functions. In the final section, this operational calculus is applied for analytical treatment of some initial-value problems for the single- and multi-term fractional differential equations with the GFDs.   

\vspace*{-3pt} 

\section{The basic rings of functions}
\label{sec:2}

\setcounter{section}{2}
\setcounter{equation}{0}\setcounter{theorem}{0}

We start with a short discussion of the ring $\mathcal{R}$ of functions introduced by Mikusi\'nski as a starting point for his operational calculus for the first order derivative. For more details see \cite{Mik59}. 

The elements of the ring $\mathcal{R}$ are the functions from the space $C[0,+\infty)$. The ring operations are the usual addition of the functions and their Laplace convolution that plays the role of multiplication:
\begin{equation}
\label{2-1}
f\, + \, g\, = \, h,\ \ h(t)=f(t)+g(t),\ \ f,g\in C[0,+\infty),
\end{equation}
\begin{equation}
\label{2-2}
f\,* \, g\, =\, h,\ \ h(t)=(f*g)(t) = \int_0^{t}\, f(\tau)g(t-\tau)\, d\tau,\ \ f,g\in C[0,+\infty).
\end{equation}
It is an easy exercise to verify that $f+g,\, f*g \in C[0,+\infty)$ if $f,g\in C[0,+\infty)$ and the triple $\mathcal{R} = (C[0,+\infty),+,*)$ satisfies all ring axioms. 

In the Mikusi\'nski operational calculus, the integration operator 
\begin{equation}
\label{2-3}
(I^1_{0+}\, f)(t)  = (\{1\}*f)(t) = \int_0^t\, f(\tau)\, d\tau,\ \  f \in C[0,+\infty)
\end{equation}
plays a very important role. Following  Mikusi\'nski, by $\{a\}$ we denote a function that is identically equal to $a$ on the whole interval $[0,+\infty)$. This function is of course not the same as the number $a,\ a\in \R$. 

According to the first fundamental theorem of calculus, the first order derivative is a left inverse operator to the first order integral, i.e.,
\begin{equation}
\label{2-4}
(D^1\, I^1_{0+}\, f)(t) = \frac{d}{dt}\, \int_0^t\, f(\tau)\, d\tau \, = \, f(t),\  f \in C[0,+\infty).
\end{equation}
The second fundamental theorem of calculus says that
\begin{equation}
\label{2-5}
(I^1_{0+}\, D^1\,  f)(t) = \int_0^t\, f^\prime(\tau)\, d\tau \, = \, f(t)-f(0),\  f \in C^1[0,+\infty).
\end{equation}
In the framework of the Mikusi\'nski operational calculus, the  formula \eqref{2-5} is employed to represent the first order derivative as multiplication in the field of convolution quotients. These matters are discussed in detail in Section \ref{sec:5} for the GFIs and the GFDs with the Sonine kernels. 

Another important ingredient of the Mikusi\'nski operational calculus is a famous theorem of Titchmarsh from \cite{Tit26} (see also \cite{Mik59} and \cite{Mik53} for other variants of its proof). This theorem ensures that the ring $\mathcal{R} = (C[0,+\infty),+,*)$ is divisors free, i.e.,  
\begin{equation}
\label{2-6}
(f*g)(t) \, \equiv \, 0,\ t\ge 0 \ \Rightarrow \ f(t) \, \equiv \, 0,\ t\ge 0 \ \mbox{ or } \ g(t) \, \equiv \, 0,\ t\ge 0. 
\end{equation}

As already mentioned in the Introduction, the Mikusi\'nski operational calculus for the first order derivative has been extended to the case of the fractional derivatives of different types. For development of these operational calculi, a larger  ring of functions compared to $\mathcal{R}$ was employed. In the rest of this Section, we provide a construction of this ring that we denote by $\mathcal{R}_{-1}$ and mention some of its properties. For the proofs of these results we refer the interested readers to \cite{LucGor99}.  

The elements of $\mathcal{R}_{-1}$ are continuous functions that can have an integrable singularity of a power function type at the point zero. We denote this space by $C_{-1}(0,+\infty)$:
\begin{equation}
\label{2-7}
f\in C_{-1}(0,+\infty)\ \ \Leftrightarrow \ \  f(t)=t^{p-1}f_1(t),\ t>0,\ p>0,\ f_1\in C[0,+\infty).
\end{equation}
A family of the spaces $C_\alpha(0,+\infty),\ \alpha \ge -1$  was first introduced in \cite{Dim66} for construction of an operational calculus for the hyper-Bessel differential operator. These spaces were then employed in the framework of the operational calculi for different fractional derivatives in \cite{BasLuc95, LucGor99, LucYak94, YakLuc94} and other related publications. 

\begin{theorem}[\cite{LucGor99}]
\label{t1}
The   triple   $\mathcal{R}_{-1} = (C_{-1}(0,+\infty),+,*)$     with  the  usual
addition $+$ and  multiplication $*$ in form of  the  Laplace convolution \eqref{2-2} is a commutative ring without divisors of zero.
\end{theorem}

Because of the inclusion $C[0,+\infty) \subset C_{-1}(0,+\infty)$ and because a sum and a convolution of any two continuous functions are again continuous functions, the ring $\mathcal{R}$ is a sub-ring of the ring $\mathcal{R}_{-1}$, i.e., $\mathcal{R} \subset \mathcal{R}_{-1}$.

In the next sections, we also need some other sub-rings of $\mathcal{R}_{-1}$. Their elements are functions  from the sub-spaces $C_{-1}^m(0,+\infty),\ m\in \N_0=\N \cup\{0\}$  of the space  $C_{-1}(0,+\infty)$: 
\begin{equation}
\label{2-7-1}
f\in C_{-1}^m(0,+\infty)\ \ \Leftrightarrow \ \  f^{(m)} \in C_{-1}(0,+\infty).
\end{equation}

The spaces $C_{-1}^m(0,+\infty)$ were introduced and studied in \cite{LucGor99}. Here we just  mention some relevant properties of these spaces, for the proofs see \cite{LucGor99}. 

\begin{theorem}[\cite{LucGor99}]
\label{t2}
For the family of spaces $C_{-1}^m(0,+\infty),\ m\in \N_0$, the following statements hold true:
\vskip 0.1cm 

\noindent
1) $C_{-1}^0(0,+\infty)\, \equiv \, C_{-1}(0,+\infty)$.
\vskip 0.1cm 

\noindent
2) $C_{-1}^m(0,+\infty),\ m\in \N_0$ is a vector  space over the field $\R$ (or $\C$).
\vskip 0.1cm  

\noindent
3) If $f\in C_{-1}^m(0,+\infty)$ with $m\ge 1$, 
then $f^{(k)}(0+) := \lim\limits_{t\to 0+} f^{(k)}(t) <+\infty,\ 0\le k\le 
m-1$, and the function
$$
\tilde f (t) = 
\begin{cases}f(t), & t>0, \\
f(0+), & t=0
\end{cases}
$$
belongs to the space $C^{m-1}[0,+\infty)$.

\vskip 0.1cm  

\noindent
4) If $f\in C_{-1}^m(0,+\infty)$ with $m\ge 1$, then $f \in 
C^{m}(0,+\infty)\cap C^{m-1}[0,+\infty)$.
\vskip 0.1cm  

\noindent
5) For any function $f\in C_{-1}^m(0,+\infty),\  m\ge 1$, the following representation holds true:
$$
f(t) = (I^m_{0+} \phi)(x) + \sum \limits_{k=0}^{m-1} f^{(k)}(0)\frac{t^k}{k!},\ t\ge 
0,\ \phi \in C_{-1}(0,+\infty).
$$
\vskip 0.1cm  

\noindent
6) Let $f\in C_{-1}^m(0,+\infty), \ m\in \N_0$,
$f(0)=\dots=f^{(m-1)}(0)=0$ and $g\in C_{-1}^1(0,+\infty)$. Then the 
Laplace convolution $h(t) = (f*g)(t)$
belongs to the  space $C_{-1}^{m+1}(0,+\infty)$ and $h(0)=\dots=h^{(m)}(0)=0$.
\end{theorem}


In the rest of this section, we provide some examples of the functions from the space $C_{-1}(0,+\infty)$ and the formulas for their multiplication in the ring $\mathcal{R}_{-1}$ (Laplace convolution).

We start with the following well-known and important formula:
\begin{equation}
\label{2-9}
(h_{\alpha} \, * \, h_\beta)(t) \, = \, h_{\alpha+\beta}(t),\ \alpha,\beta >0,\ t>0,
\end{equation}
where the power function $h_\beta$ is defined by
\begin{equation}
\label{h}
h_{\beta}(t) := \frac{t^{\beta -1}}{\Gamma(\beta)},\ \beta >0.
\end{equation}
For $\alpha,\beta >0$, the inclusions  $h_{\alpha},\ h_{\beta}, \, h_{\alpha+\beta} \in C_{-1}(0,+\infty)$ hold true. The formula \eqref{2-9} is an easy consequence from the well-known representation of the Euler Beta-function in terms of the Gamma-functions. 

The same representation  ensures a more general relation 
\begin{equation}
\label{2-10}
(h_{\alpha,\rho} \, * \, h_{\beta,\rho})(t) \, = \, h_{\alpha+\beta,\rho}(t),\ \alpha,\beta >0,\ \rho \in \R,\  t>0
\end{equation}
with $h_{\alpha,\rho} \in C_{-1}(0,+\infty)$ defined as follows:
\begin{equation}
\label{2-11}
h_{\alpha,\rho}(t) = \frac{t^{\alpha -1}}{\Gamma(\alpha)}\, e^{-\rho t},\ \ \alpha >0,\ \rho \in \R,\ t>0.
\end{equation}

Because $\mathcal{R}_{-1}$ is a ring, the integer order convolution powers of a function $\kappa \in C_{-1}(0,+\infty)$ are well defined:
\begin{equation}
\label{2-12}
\kappa^n :=\underbrace{\kappa*\ldots\ * \kappa}_n,\ n\in \N.
\end{equation}
For the power function $h_\alpha$,  the formula \eqref{2-9} implicates the representation
\begin{equation}
\label{2-13}
h_{\alpha}^n(t) = h_{n\alpha}(t),\ n\in \N.
\end{equation}
Its important and well-known particular case is given by the relation
\begin{equation}
\label{2-13-1}
\{1\}^n(t) = h_{1}^n(t)= h_{n}(t)=\frac{t^{n-1}}{\Gamma(n)} = \frac{t^{n-1}}{(n-1)!},\ n\in \N.
\end{equation}

On the ring $\mathcal{R}_{-1}$, an analogy of the binomial formula in form
\begin{equation}
\label{2-14}
(f\, +\, g)^n = \sum_{i=0}^n \binom{n}{i}\, f^i\, *\, g^{n-i}
\end{equation}
evidently holds true. This formula can be used to get the series representations of the integer order convolution powers of the sums of two or several functions from  the space $C_{-1}(0,+\infty)$. 

As an example, let us consider a sum of two power functions
\begin{equation}
\label{2-15}
\kappa(t) = h_{1-\beta+\alpha}(t)\, +\, h_{1-\beta}(t),\ 0<\alpha < \beta <1.
\end{equation}
Then $\kappa \in C_{-1}(0,+\infty)$ and we get the following chain of equations:
$$
\kappa^n(t) = (h_{1-\beta+\alpha}(t)\, +\, h_{1-\beta}(t))^n = \sum_{i=0}^n \binom{n}{i}\, h_{1-\beta+\alpha}^i(t)\, *\, h_{1-\beta}^{n-i}(t) = 
$$
\begin{equation}
\label{2-16}
\sum_{i=0}^n \binom{n}{i} (h_{i(1-\beta+\alpha)}*h_{(n-i)(1-\beta)})(t) = 
\sum_{i=0}^n \binom{n}{i} h_{i\alpha + n(1-\beta)}(t). 
\end{equation}




\vspace*{-3pt} 

\section{The GFIs and GFDs with the Sonine kernels}
\label{sec:4}

\setcounter{section}{3}
\setcounter{equation}{0}\setcounter{theorem}{0}

In this section, we address the GFIs and the GFDs  with the special Sonine kernels on the spaces of functions $C_{-1}^m(0,+\infty), \ m\in \N_0$ discussed in Section \ref{sec:2} and on their sub-spaces. 

First we shortly discuss the Sonine kernels. For more details we refer the readers to \cite{Luc21,Sam}. In \cite{Son}, Sonine extended  Abel's method for analytical treatment of the Abel integral equation to more general integral equations of convolution type. He generalized the formula
\begin{equation}
\label{3-1}
(h_{\alpha}\, * \, h_{1-\alpha})(t) = h_1(t) = \{1 \},\ 0<\alpha<1,\  t>0
\end{equation}
that is in fact the most essential ingredient of Abel's solution to the arbitrary kernels $\kappa$ and $k$ that satisfy the relation
\begin{equation}
\label{3-2}
(\kappa \, *\, k )(t) = \{1 \},\  t>0.
\end{equation}
Nowadays, such functions are called the Sonine kernels and the relation \eqref{3-2} is called the Sonine condition. For a Sonine kernel $\kappa$, the kernel $k$ that satisfies the Sonine condition \eqref{3-2} is called an associate kernel. Vise versa, $\kappa$ is then an associated kernel to $k$. The set of the Sonine kernels will be denoted by $\mathcal{S}$.

In \cite{Han20}, the Sonine condition was shown to be sufficient for validity of the 1st fundamental theorem of FC for the GFD and the corresponding GFI with the kernels $k$ and $\kappa$, respectively. As mentioned in \cite{Sam} (see also \cite{Han20}), any locally integrable Sonine kernel possesses a singularity at the point zero.  On the other hand, the kernels of any fractional integral and derivative should be singular (see \cite{DGGS} for a detailed discussion of this topic). Thus, the GFIs and GFDs with the Sonine kernels are worth to be investigated. 

In the rest of this paper, we deal with the GFIs and GFDs with the Sonine kernels that have an integrable singularity of the power function type at the point zero.
 
\begin{definition}
\label{dd2}
Let $\kappa,\, k $ be a pair of the Sonine kernels that belong to the space
\begin{equation}
\label{subspace}
 C_{-1,0}(0,+\infty) \, = \, \{f:\ f(t) = t^{p-1}f_1(t),\ t>0,\ 0<p<1,\ f_1\in C[0,+\infty)\}.
\end{equation}

The set of such Sonine kernels will be denoted by $\mathcal{S}_{-1}$:
\begin{equation}
\label{Son}
(\kappa,\, k \in \mathcal{S}_{-1} ) \ \Leftrightarrow \ (\kappa,\, k \in C_{-1,0}(0,+\infty))\wedge ((\kappa\, *\, k)(t) \, = \, \{1\}).
\end{equation}
\end{definition}

An immediate and important consequence from Theorem \ref{t1} (the ring $\mathcal{R}_{-1}$ is divisors free) is uniqueness of an associate kernel to any  Sonine kernel from $\mathcal{S}_{-1}$. It is also worth mentioning that the kernels of the most time-fractional derivatives and integrals introduced so far are from the set $\mathcal{S}_{-1}$ and thus the constructions presented in the next sections are applicable both for the known and for many new FC operators.  In particular, the kernels of the Riemann-Liouville fractional integral and derivative, $h_{\alpha}$ and $h_{1-\alpha}$, respectively, belong to $\mathcal{S}_{-1}$ if $0<\alpha < 1$. 

The general classes of the Sonine kernels pairs from $\mathcal{S}_{-1}$ are given by the formulas
\begin{equation}
\label{3-3}
\kappa(t) = h_{\alpha}(t) \cdot \, \kappa_1(t),\ \kappa_1(t)=\sum_{k=0}^{+\infty}\, a_k t^k, \ a_0 \not = 0,\ 0<\alpha <1,
\end{equation}
where $\kappa_1=\kappa_1(t)$ has the infinite (\cite{Son}) or finite convergence radius $r>0$ (\cite{Wick}) and
\begin{equation}
\label{3-4}
k(t) = h_{1-\alpha}(t) \cdot k_1(t),\ k_1(t)=\sum_{k=0}^{+\infty}\, b_k t^k, 
\end{equation}
where the coefficients $b_k,\ k\in \N_0$ are uniquely determined by the coefficients $a_k,\ k\in \N_0$ as solutions to the following triangular system of linear equations:
\begin{equation}
\label{3-5}
a_0b_0 = 1,\ \sum_{k=0}^n\Gamma(k+1-\alpha)\Gamma(\alpha+n-k)a_{n-k}b_k = 0,\ n\ge 1.
\end{equation}
A famous particular case of the kernels in form \eqref{3-3}, \eqref{3-4} was deduced by Sonine in \cite{Son}:
\begin{equation}
\label{Bess}
\kappa(t) = (\sqrt{t})^{\alpha-1}J_{\alpha-1}(2\sqrt{t}),\ 
k(t) = (\sqrt{t})^{-\alpha}I_{-\alpha}(2\sqrt{t}),\ 0<\alpha <1,
\end{equation}
where 
$$
J_\nu (t) = \sum_{k=0}^{+\infty} \frac{(-1)^k(t/2)^{2k+\nu}}{k!\Gamma(k+\nu+1)},\ 
I_\nu (t) = \sum_{k=0}^{+\infty} \frac{(t/2)^{2k+\nu}}{k!\Gamma(k+\nu+1)}
$$
are the Bessel and the modified Bessel functions, respectively. Another example of this type is the following pair of the Sonine kernels 
(\cite{Zac08}):
\begin{equation}
\label{3-6} 
\kappa(t) = 
h_{\alpha,\rho}(t),\ 0<\alpha <1,\ \rho\ge 0,
\end{equation}
\begin{equation}
\label{3-7} 
k(t) = h_{1-\alpha,\rho}(t) \, +\, \rho\, \int_0^t h_{1-\alpha,\rho}(\tau)\, d\tau,
\end{equation}
where the function $h_{\alpha,\rho}$ is defined by the formula \eqref{2-11}. 
The condition $\rho \ge 0$ allows us to represent the kernel \eqref{3-7}  as follows:
\begin{equation}
\label{3-7-1} 
k(t) = h_{1-\alpha,\rho}(t) \, +\, \frac{\rho^\alpha}{\Gamma(1-\alpha)}\, \gamma(1-\alpha,\rho t),
\end{equation}
where $\gamma$ stands for the incomplete Gamma-function:
$$
\gamma(\beta,t) = \int_0^t \tau^{\beta -1} e^{-\tau}\, d\tau,\ t>0.
$$

In \cite{Koch11}, an important special class of the Sonine kernels was introduced. This class consists of all functions $k = k(t),\ t>0$ that satisfy the conditions K1)-K4) listed in the Introduction. As already mentioned, the solutions to the fractional differential equations with the GFDs of type \eqref{FDC} with the kernels from this class  behave as the ones of the evolution equations  (\cite{Koch11}). 

In \cite{Koch11}, the Sonine kernels were treated in terms of their Laplace transforms and using the Sonine condition \eqref{3-2} in the Laplace domain: 
\begin{equation}
\label{Laplace-Sonine} 
\tilde \kappa(p) \cdot \tilde k(p)  = \frac{1}{p},\ \Re(p)>p_{\kappa,k} \in \R,
\end{equation}
provided the Laplace transforms $\tilde \kappa,\ \tilde k$ of the functions $\kappa$ and $k$ do exist.

In \cite{Han20}, another class of the Sonine kernels was introduced in terms of the completely monotone functions. According to \cite{Han20}, any singular (unbounded in a neighborhood of the point zero) locally integrable
completely monotone function is a Sonine kernel. A typical example of such Sonine kernels is as follows (\cite{Han20}):
\begin{equation}
\label{3-8} 
\kappa(t) = h_{1-\beta+\alpha}(t)\, +\, h_{1-\beta}(t),\ 0<\alpha < \beta <1,
\end{equation}
\begin{equation}
\label{3-9} 
k(t) = t^{\beta -1}\, E_{\alpha,\beta}(-t^\alpha),
\end{equation}
where $E_{\alpha,\beta}$ stands for the two-parameters Mittag-Leffler function defined by the following convergent series:
\begin{equation}
\label{ML}
E_{\alpha,\beta}(z) = \sum_{k=0}^{+\infty} \frac{z^k}{\Gamma(\alpha\, k + \beta)},\ \alpha >0,\ \beta,z\in \C.
\end{equation}

In the rest of this section, following \cite{Luc21}, we discuss some important properties of the GFIs and GFDs with the Sonine kernels from $\mathcal{S}_{-1}$.

\begin{definition}
\label{d-GFI}
Let $\kappa \in \mathcal{S}_{-1}$. The GFI with the kernel $\kappa$ is defined by the formula
\begin{equation}
\label{GFI}
(\I_{(\kappa)}\, f)(t) := (\kappa\, *\, f)(t) = \int_0^t \kappa(t-\tau)f(\tau)\, d\tau,\ t>0.
\end{equation}
\end{definition}

For $\kappa(t) = h_\alpha(t),\ 0<\alpha<1$, the GFI \eqref{GFI} is the Riemann-Liouville fractional integral:
\begin{equation}
\label{GFI_1}
(\I_{(\kappa)}\, f)(t) = (I^\alpha_{0+}\, f)(t),\ t>0,
\end{equation} 
where 
\begin{equation}
\label{RLI}
\left(I^\alpha_{0+} f \right) (t)=(h_\alpha \, * \, f)(x) = \frac{1}{\Gamma (\alpha )}\int\limits_0^t(t -\tau)^{\alpha -1}f(\tau)\,d\tau,\ t>0.
\end{equation}
Taking other known Sonine kernels, other particular cases of the GFI \eqref{GFI} can be constructed. Here we mention just a few of them:
\begin{equation}
\label{GFI_2}
(\I_{(\kappa)}\, f)(t) = (I^{1-\beta+\alpha}_{0+}\, f)(t) + (I^{1-\beta}_{0+}\, f)(t)  ,\ t>0
\end{equation} 
with the Sonine kernel $\kappa$ defined by \eqref{3-8},
\begin{equation}
\label{GFI_3}
(\I_{(\kappa)}\, f)(t) = \frac{1}{\Gamma(\alpha)}\int_0^t (t-\tau)^{\alpha -1}\, e^{-\rho (t-\tau)}\, f(\tau)\, d\tau,\ t>0
\end{equation} 
with the Sonine kernel $\kappa$ defined by \eqref{3-6}, and 
\begin{equation}
\label{GFI_4}
(\I_{(\kappa)}\, f)(t) = \int_0^t (\sqrt{t-\tau})^{\alpha-1}J_{\alpha-1}(2\sqrt{t-\tau})f(\tau)\, d\tau,\ t>0
\end{equation}
with the Sonine kernel $\kappa$ defined by \eqref{Bess}.

The properties of the GFI \eqref{GFI} on the space $C_{-1}(0,+\infty)$ result from the known properties of the Laplace convolution and Theorem \ref{t2}. In particular, the following relations hold true on the space $C_{-1}(0,+\infty)$  (\cite{Luc21}): 
\begin{equation}
\label{GFI-map}
\I_{(\kappa)}:\, C_{-1}(0,+\infty)\, \rightarrow C_{-1}(0,+\infty),
\end{equation}
\begin{equation}
\label{GFI-com}
\I_{(\kappa_1)}\, \I_{(\kappa_2)} = \I_{(\kappa_2)}\, \I_{(\kappa_1)},\ \kappa_1,\, \kappa_2 \in \mathcal{S}_{-1},
\end{equation}
\begin{equation}
\label{GFI-index}
\I_{(\kappa_1)}\, \I_{(\kappa_2)} = \I_{(\kappa_1*\kappa_2)},\  \kappa_1,\, \kappa_2 \in \mathcal{S}_{-1}.
\end{equation}

\begin{definition}
\label{d-GFD}
The GFDs of the Riemann-Liouville and Caputo types that correspond to the GFI \eqref{GFI} with a Sonine kernel $\kappa \in \mathcal{S}_{-1}$ are defined as follows: 
\begin{equation}
\label{GFDL}
(\D_{(k)}\, f) (t) = \frac{d}{dt}(k\, * \, f)(t),\ t>0,
\end{equation}
\begin{equation}
\label{GFDC}
( _*\D_{(k)}\, f) (t) =  (\D_{(k)}\, f) (t) - f(0)k(t),\ t>0,
\end{equation}
where the kernel $k \in \mathcal{S}_{-1}$ is the Sonine kernel associate to the kernel $\kappa$. 
\end{definition}

\begin{remark}
\label{r1}
It is worth mentioning that both the integral operators in form \eqref{GFI} and the integro-differential operators in form \eqref{GFDL} and \eqref{GFDC} have been already introduced and studied in the literature on different spaces of functions and with different kernels. The novelty of our approach is in treating these operators with the kernels from a very large class $\mathcal{S}_{-1}$ on the suitable spaces of functions $C_{-1}^m(0,+\infty),\ m\in \N_0$. As shown in the rest of this paper, the integral operators \eqref{GFI} and the integro-differential operators  \eqref{GFDL} and \eqref{GFDC} with the kernels from $\mathcal{S}_{-1}$ still possess the main characteristics that justify calling them the GFIs and the GFDs, respectively.
\end{remark}

Simple calculations lead to the following interpretation of  the GFD \eqref{GFDC} in the Caputo sense  as a regularized GFD \eqref{GFDL} in the Riemann-Liouville sense (see Definition 3.2  of the Caputo fractional derivative in \cite{Diet}):
\begin{equation}
\label{GFDC_new}
( _*\D_{(k)}\, f) (t) =  (\D_{(k)}\, [f(\cdot)-f(0)]) (t),\ t>0.
\end{equation}

In this paper, we mainly deal with the GFDs in form \eqref{GFDC} on the spaces $C_{-1}^m(0,+\infty),$ $m\in \N$ and their sub-spaces. For the functions from $C_{-1}^1(0,+\infty)$, Theorem \ref{t2} and the known formula for differentiation of the integrals depending on parameters lead to an important representation of the GFD \eqref{GFDL} in the Riemann-Liouville sense:
\begin{equation}
\label{GFDL-1}
(\D_{(k)}\, f) (t) = (k\, * \, f^\prime)(t) + f(0)k(t),\ t>0.
\end{equation}
A a consequence, for $f\in C_{-1}^1(0,+\infty)$, the GFD \eqref{GFDC} in the Caputo sense can be rewritten in the form
\begin{equation}
\label{GFDC_1}
( _*\D_{(k)}\, f) (t) =  (k\, * \, f^\prime)(t) = (\I_{(k)}\, f^\prime)(t),\ t>0. 
\end{equation}

For the Riemann-Liouville fractional integral \eqref{RLI} with the kernel $\kappa(t) = h_{\alpha}(t)$, the associate Sonine kernel is $k(t)=h_{1-\alpha}(t)$ and thus the GFDs \eqref{GFDL} and \eqref{GFDC} are the Riemann-Liouville and the Caputo  fractional derivatives of order $\alpha,\ 0< \alpha < 1$:
\begin{equation}
\label{RLD}
(\D_{(k)}\, f) (t) = \frac{d}{dt}(h_{1-\alpha}\, * \, f)(t)  = (D^\alpha\, f)(t),\ t>0,
\end{equation}
\begin{equation}
\label{CD}
( _*\D_{(k)}\, f) (t) = \frac{d}{dt}(h_{1-\alpha}\, * \, f)(t) - f(0)h_{1-\alpha}(t) = ( _*D^\alpha\, f)(t),\ t>0. 
\end{equation}

Other conventional particular cases of the GFDs \eqref{GFDL} and \eqref{GFDC} are the multi-term fractional derivatives and the fractional derivatives of the
distributed order. They are the operators \eqref{GFDL} and \eqref{GFDC} with the kernels 
\begin{equation}
\label{multi}
k(t) = \sum_{k=1}^n a_k\, h_{1-\alpha_k}(t),
\ \  0 < \alpha_1 <\dots < \alpha_n < 1,\ a_k\in \R,\ k=1,\dots,n
\end{equation}
and
\begin{equation}
\label{distr}
k(t) = \int_0^1 h_{1-\alpha}(t)\, d\rho(\alpha),
\end{equation}
respectively, where $\rho$ is a Borel measure defined on the interval $[0,\, 1]$.

Of course, one can construct many other GFDs based on the known Sonine kernels from $\mathcal{S}_{-1}$. In the next sections, we employ some particular cases of the GFD in the Caputo sense that correspond to the  GFIs \eqref{GFI_1}, \eqref{GFI_2}-\eqref{GFI_4}. To make the formulas more compact, we present them in form  \eqref{GFDC_1} (of course, they also  can be represented in form \eqref{GFDC}). 

The GFD that corresponds to the GFI \eqref{GFI_2} has the Mittag-Leffler function \eqref{3-9} in the kernel:
\begin{equation}
\label{GFD_2}
( _*\D_{(k)}\, f)(t) = \int_0^t \tau^{\beta -1}\, E_{\alpha,\beta}(-\tau^\alpha)\,   f^\prime(t-\tau)\, d\tau,\ 0<\alpha < \beta <1,\ t>0.
\end{equation} 

The Sonine kernel associate to the  kernel $\kappa(t) = h_{\alpha,\rho}(t),\ 0<\alpha<1,\ \rho\ge 0$ is given by the formula \eqref{3-7-1}. 
Thus, the GFD associate to the GFI \eqref{GFI_3} has the form 
\begin{equation}
\label{GFD_3}
( _*\D_{(k)}\, f)(t) = \int_0^t \left( h_{1-\alpha,\rho} (\tau) \, +\, \frac{\rho^\alpha}{\Gamma(1-\alpha)}\, \gamma(1-\alpha,\rho \tau) \right)\, f^\prime(t-\tau)\, d\tau,\ t>0.
\end{equation} 
Finally, the GFD that corresponds to the GFI \eqref{GFI_4} takes the form (see the Sonine pair \eqref{Bess}): 
\begin{equation}
\label{GFD_4}
( _*\D_{(k)}\, f)(t)  = \int_0^t 
(\sqrt{\tau})^{-\alpha}I_{-\alpha}(2\sqrt{\tau})f^\prime(t-\tau)\, d\tau,\ t>0.
\end{equation}

In the rest of this section, we present some results that justify referring to the integral operator \eqref{GFI} and the integro-differential operators \eqref{GFDL} and \eqref{GFDC} as to the GFI and the GFDs, respectively.  For the proofs of these results we refer to \cite{Luc21}. 

\begin{theorem}[1st fundamental theorem of FC for the GFD]
\label{t3}
Let $k \in \mathcal{S}_{-1}$ be  the Sonine kernel associate to the kernel $\kappa$.  

Then the GFD  \eqref{GFDL} is a left inverse operator to the GFI \eqref{GFI} on the space $C_{-1}(0,+\infty)$: 
\begin{equation}
\label{FTL}
(\D_{(k)}\, \I_{(\kappa)}\, f) (t) = f(t),\ f\in C_{-1}(0,+\infty),\ t>0,
\end{equation}
and the GFD \eqref{GFDC} is a left inverse operator to the GFI \eqref{GFI} on the space $C_{-1,k}(0,+\infty)$: 
\begin{equation}
\label{FTC}
( _*\D_{(k)}\, \I_{(\kappa)}\, f) (t) = f(t),\ f\in C_{-1,k}(0,+\infty),\ t>0,
\end{equation}
where
$$
C_{-1,k}(0,+\infty):= \{f:\ f(t)=(\I_{(k)}\, \phi)(t),\ \phi\in C_{-1}(0,+\infty)\}.
$$ 
\end{theorem}

As mentioned in \cite{Luc21}, the space $C_{-1,k}(0,+\infty)$ can be also characterized as follows:
$$
C_{-1,k}(0,+\infty) = 
 \{f:\ \I_{(\kappa)} f \in C_{-1}^1(0,+\infty)\, \wedge \, (\I_{(\kappa)}\, f)(0)  = 0\}.
$$
\begin{theorem}[2nd fundamental theorem of FC for the GFD]
\label{t4}
Let $k \in \mathcal{S}_{-1}$ be  the Sonine kernel associate to the kernel $\kappa$ and $f\in C_{-1}^1(0,+\infty)$.   

Then the relations
\begin{equation}
\label{2FTC}
(\I_{(\kappa)}\,  _*\D_{(k)}\, f) (t) = f(t)-f(0),\ t>0,
\end{equation}
\begin{equation}
\label{2FTL}
(\I_{(\kappa)}\, \D_{(k)}\,  f) (t) = f(t),\ t>0
\end{equation}
hold true.
\end{theorem}

Having in mind development of an operational calculus for the GFD \eqref{GFDC}, we consider the compositions of the GFIs and construct the suitable fractional derivatives. 

\begin{definition}
\label{d1}
Let $\kappa \in \mathcal{S}_{-1}$.  The $n$-fold GFI ($n \in \N$) is defined as a composition of $n$ GFIs  with the kernel $\kappa$:
\begin{equation}
\label{GFIn}
(\I_{(\kappa)}^n\, f)(t) := (\underbrace{\I_{(\kappa)} \ldots \I_{(\kappa)}}_n\, f)(x) = (\kappa^n\, *\, f)(t),\ \kappa^n:= \underbrace{\kappa* \ldots * \kappa}_n,\ t>0.
\end{equation}
\end{definition}

Evidently, the kernel $\kappa^n,\ n\in \N$ is from the space 
$C_{-1}(0,+\infty)$, but it is not always a Sonine kernel. For instance, the $n$-fold GFI \eqref{GFIn} with the kernel $\kappa(t) = h_\alpha(t)$ is the Riemann-Liouville fractional integral of the order $n\alpha,\ n\in \N$. Its kernel $h_{n\alpha}$ is a Sonine kernel only under the condition 
$\alpha<1/n$. However, the Riemann-Liouville fractional integral is interpreted as a fractional integral for any $\alpha >0$, i.e., also in the case when its kernel  is not singular at the point zero and thus not a Sonine kernel. This is justified by existence of an appropriate fractional derivative, the Riemann-Liouville fractional derivative, that is connected with the Riemann-Liouville fractional integral through the 1st and the 2nd fundamental theorems of FC. In the rest of this section, we introduce an analogous construction for the GFI \eqref{GFI} that justifies calling the operator \eqref{GFIn} a fractional integral. 

\begin{definition}
\label{d2}
Let $k$  be a Sonine kernel from  $\mathcal{S}_{-1}$. The $n$-fold GFD ($n \in \N$) in the Riemann-Liouville sense is defined as follows:
\begin{equation}
\label{GFDLn}
(\D_{(k)}^n\, f)(t) := \frac{d^n}{dt^n} ( k^n * f)(t),\, k^n:= \underbrace{k* \ldots * k}_n,\ t>0.
\end{equation}
\end{definition}

The $n$-fold GFD is a generalization of the Riemann-Liouville fractional derivative of an arbitrary order $\alpha > 0$ to the case of the general Sonine kernels. For instance, for the Sonine kernels pair $\kappa(t) = h_{\alpha}(t),\ k(t) = h_{1-\alpha}(t),\ 0<\alpha < 1$, the two-fold GFD \eqref{GFDLn} has the kernel $k(t) = h_{2-2\alpha}(t)$ and thus it can be represented in terms of the Riemann-Liouville fractional derivative:
\begin{equation}
\label{RLn}
(\D_{(k)}^2\, f)(t) = (D_{0+}^{2\alpha}\, f)(t) = 
\begin{cases} 
\frac{d^2}{dt^2}(I^{2-2\alpha}_{0+}\, f)(t),& \frac{1}{2} < \alpha <1,\ t>0, \\
 \frac{d}{dt}(I^{1-2\alpha}_{0+}\, f)(t),& 0 < \alpha \le \frac{1}{2},\ t>0. 
 \end{cases}
\end{equation}

The Definition \ref{d2} immediately implicates  the following recurrent formula for the $n$-fold GFD \eqref{GFDL}:
\begin{equation}
\label{GFDLn_rec}
(\D_{(k)}^n\, f) (t) = 
\frac{d}{dt}\frac{d^{n-1}}{dt^{n-1}} \left( (k^{n-1}*(k * f))(t)\right) = 
\frac{d}{dt} (\D_{(k)}^{n-1}\, (k*f)) (t),\ n\in \N.
\end{equation}
This formula can be employed for extension of the results known for the GFD \eqref{GFDL} to the case of the $n$-fold GFD \eqref{GFDLn}. Let us note that the formula \eqref{GFDLn_rec} is valid also for $n=1$ provided we define $(\D_{(k)}^{0}\, f)(t) := f(t),\ t>0$. 

In particular, the following important generalization of the representation \eqref{GFDL-1} of the GFD \eqref{GFDL} to the case of the $n$-fold GFD \eqref{GFDLn} has been deduced: 

\begin{theorem}[\cite{Luc21}]
\label{t5}
Let $k$ be a Sonine kernel from $\mathcal{S}_{-1}$. 

Then the representation 
\begin{equation}
\label{GFDLn-1}
(\D_{(k)}^n\, f) (t) = (k^n\, * \, f^{(n)})(t) + \sum_{j=0}^{n-1}f^{(j)}(0) 
\frac{d^{n-j-1}}{dt^{n-j-1}} k^n(t),\ t>0
\end{equation}
holds true for any function $f\in C_{-1}^n(0,+\infty)$. 
\end{theorem}

It is worth mentioning that the formula \eqref{GFDLn-1} with the kernel $k(t) = h_{1-\alpha}(t),\ 0<\alpha <1$ connects the Riemann-Liouville and the Caputo fractional derivatives of the order $n\alpha$ (see Lemma 3.4 from \cite{Diet}):
\begin{equation}
\label{GFDLn-1-RL}
(D_{0+}^{n\alpha}\, f) (t) = ( _*D^{n\alpha}\, f) (t) + \sum_{j=0}^{n-1}f^{(j)}(0) 
h_{j+1-n\alpha}(t),\ t>0.
\end{equation}

As in the case of the GFD \eqref{GFDL}, the representation 
\eqref{GFDLn-1} is used to define a Caputo type $n$-fold GFD. 

\begin{definition}
\label{d3}
Let $k$  be a Sonine kernel from  $\mathcal{S}_{-1}$. The $n$-fold GFD ($n \in \N$) in the Caputo sense is defined as follows:
\begin{equation}
\label{GFDCn}
( _*\D_{(k)}^n\, f)(t) := (\D_{(k)}^n\, f) (t) - \sum_{j=0}^{n-1}f^{(j)}(0) 
\frac{d^{n-j-1}}{dt^{n-j-1}} k^n(t),\ t>0.
\end{equation}
\end{definition}
An equivalent definition of the $n$-fold Caputo type GFD can be given in form of a regularized $n$-fold GFD in the Riemann-Liouville sense (see Definition 3.2 of the Caputo fractional derivative in \cite{Diet}):
\begin{equation}
\label{GFDCn-new}
( _*\D_{(k)}^n\, f)(t) = \left(\D_{(k)}^n\, \left[f(\cdot) - \sum_{j=0}^{n-1}f^{(j)}(0)\, h_{j+1}(\cdot) \right]\right)(t),\ t>0.
\end{equation}
Indeed, for $k\in \mathcal{S}_{-1}$, we have the following chain of equations:
$$
\left(\D_{(k)}^n\, \left[f(\cdot) - \sum_{j=0}^{n-1}f^{(j)}(0)\, h_{j+1}(\cdot)\right]\right)(t) = 
$$
$$
(\D_{(k)}^n\, f)(t) - \sum_{j=0}^{n-1} f^{(j)}(0)(\D_{(k)}^n\, h_{j+1})(t) =
$$
$$
(\D_{(k)}^n\, f)(t) - \sum_{j=0}^{n-1} f^{(j)}(0)\, \frac{d^n}{dt^n} (k^n\, *\, h_{j+1})(t) = 
$$
$$
(\D_{(k)}^n\, f)(t) - \sum_{j=0}^{n-1} f^{(j)}(0)\, \frac{d^n}{dt^n} (I_{0+}^{j+1}\, k^n)(t) =
$$
$$
(\D_{(k)}^n\, f) (t) - \sum_{j=0}^{n-1}f^{(j)}(0) 
\frac{d^{n-j-1}}{dt^{n-j-1}} k^n(t),\ t>0.
$$

According to Theorem \ref{t5}, on the space $C_{-1}^n(0,+\infty)$,  the $n$-fold GFD $ _*\D_{(k)}^n$ can be also represented in the form
\begin{equation}
\label{GFDCn_1}
( _*\D_{(k)}^n\, f)(t) = (k^n\, * \, f^{(n)})(t),\ t>0.
\end{equation}
As in the case of the Caputo fractional derivative \eqref{CD}, the $n$-fold GFD in form \eqref{GFDCn} makes sense also for the functions that are not $n$-times differentiable (for instance, for the functions from the space $AC^{n-1}[0,+\infty)$). However, the representation of type \eqref{GFDCn_1} for the Caputo fractional derivative is very convenient and thus widely used. 

The denotation "$n$-fold GFD" is justified by the fact that the $n$-fold GFDs in the Riemann-Liouville and in the Caputo sense  are left inverse operators to the $n$-fold GFI \eqref{GFIn}. 

\begin{theorem}[1st fundamental theorem of FC for the $n$-fold GFD]
\label{t6}
Let $k \in \mathcal{S}_{-1}$ be  the Sonine kernel associate to the kernel $\kappa$.  

Then the $n$-fold GFD  \eqref{GFDLn} is a left inverse operator to the $n$-fold GFI \eqref{GFIn} on the space $C_{-1}(0,+\infty)$: 
\begin{equation}
\label{FTLn}
(\D_{(k)}^n\, \I_{(\kappa)}^n\, f) (t) = f(t),\ f\in C_{-1}(0,+\infty),\ t>0
\end{equation}
and the GFD \eqref{GFDCn} is a left inverse operator to the  $n$-fold GFI \eqref{GFIn} on the space $C_{-1,k}^n(0,+\infty) := \{f:\ f(t)=(\I_{(k)}^n\, \phi)(t),\ \phi\in C_{-1}(0,+\infty)\}$: 
\begin{equation}
\label{FTCn}
( _*\D_{(k)}^n\, \I_{(\kappa)}^n\, f) (t) = f(t),\ f\in C_{-1,k}^n(0,+\infty),\ t>0.
\end{equation}
\end{theorem}
 
Finally, we formulate the 2nd fundamental theorem of FC for the $n$-fold GFD in the  Caputo sense (for the proof see \cite{Luc21}).

\begin{theorem}[2nd fundamental theorem of FC for the $n$-fold GFD]
\label{t7}
Let $k \in \mathcal{S}_{-1}$ be  the Sonine kernel associate to the kernel $\kappa$ and $f\in C_{-1}^n(0,+\infty)$.   

Then the relation
\begin{equation}
\label{2FTCn}
(\I_{(\kappa)}^n\,  _*\D_{(k)}^n\, f) (t) = f(t)-\sum_{j=0}^{n-1} \, f^{(j)}(0)\, h_{j+1}(t),\ t>0
\end{equation}
holds true.
\end{theorem}

\vspace*{-3pt} 

\section{Operational calculus for the GFD with a Sonine kernel}
\label{sec:5}

\setcounter{section}{4}
\setcounter{equation}{0}\setcounter{theorem}{0}

In this section, we present the basic constructions of a Mikusi\'nski type operational calculus for the GFD \eqref{GFDC} in the Caputo sense with a Sonine kernel $k\in \mathcal{S}_{-1}$. 

The GFI \eqref{GFI} acting on a function $f\in C_{-1}(0,+\infty)$ is defined 
as  multiplication of the function $f$ with the Sonine kernel $\kappa \in 
\mathcal{S}_{-1}$ on the ring $\mathcal{R}_{-1}$. 
Because the  GFD \eqref{GFDC} is a left inverse operator to the GFI \eqref
{GFI}  (Theorem \ref{t3}), a representation of the GFD as multiplication on the ring $\mathcal{R}_{-1}$ is not possible (the ring $\mathcal{R}_{-1}$ does not possess a unity element with respect to multiplication). In this section, we extend the ring $\mathcal{R}_{-1}$  to a field of convolution quotients and introduce an algebraic GFD as a multiplication on this field. 

In what follows, we closely follow the lines of development of the operational calculi for the fractional derivatives of different types (see, e.g., \cite{BasLuc95,HadLuc,LucGor99,LucSri95,LucYak94,YakLuc94}). First, we introduce a 
set of functions pairs from $C_{-1}(0,+\infty)$ 
$$
C_{-1}^2:= C_{-1}(0,+\infty) \times  (C_{-1}(0,+\infty) \setminus \{0\})
$$
and an equivalence
relation defined on this set as follows:
$$
(f_1,\, g_1)\sim     (f_2,\, g_2)    \Leftrightarrow     (f_1\, * \, 
g_2)(t) = (f_2\, * \,  g_1)(t),\ (f_1,\, g_1),\ (f_2,\, g_2) \in C_{-1}^2.
$$
The elements of the convolution quotients field are classes of equivalences $C_{-1}^2/\sim$ that we denote as quotients:
$$
\frac{f}{g}:= \{ (f_1,\, g_1) \in C_{-1}^2:\ (f_1,\, g_1) \sim (f,\, g) \}.
$$
On the set $C_{-1}^2/\sim$, the operations of addition and multiplication are defined as usual:
$$
\frac{f_1}{g_1} + \frac{f_2}{g_2}:= \frac{ f_1\, *\,  g_2\, +\, f_2\, * g_1}
{g_1\, *\, g_2},
$$
$$
\frac{f_1}{g_1}\cdot  \frac{f_2}{g_2} :=\frac{f_1\, *\,  f_2}{g_1\, *\,  g_2}.
$$
It is easy to verify that the operations $+$ and $\cdot$ introduced above do not depend on the representatives of the  equivalence classes and thus are correctly defined. 

Moreover, Theorem \ref{t1} immediately implicates the following important statement:

\begin{theorem}[\cite{LucGor99}]
\label{t8}
The triple $\mathcal{F}_{-1} = (C_{-1}^2/\sim,\ +,\ \cdot)$ is a  field that is called the field of convolution quotients.
\end{theorem}

The ring $\mathcal{R}_{-1}$ can be embedded into the field
$\mathcal{F}_{-1}$, say,  by the mapping:
\begin{equation}
\label{emb}
f  \mapsto \frac{f\, *\, \kappa}{\kappa},
\end{equation}
where $\kappa \in \mathcal{S}_{-1}$ is the kernel of the GFI \eqref{GFI}. Of course, instead of $\kappa$, any other element of $C_{-1}(0,+\infty)$ not identically equal to zero can be taken. 

On the space $C_{-1}^2/\sim$, multiplication with a scalar $\lambda \in \R$ 
(or $\lambda \in \C$) can be introduced as follows:
$$
\lambda\, \frac{f}{g}:=\frac{\lambda\, f}{g},\ \frac{f}{g}\in C_{-1}^2/\sim.
$$
Because the space of functions $C_{-1}(0,+\infty)$ is a vector space, we easily verify
that the set $C_{-1}^2/\sim$ is also a vector space. It is worth mentioning that  the constant function $\{ \lambda \}$ ($f(t)\equiv \lambda,\ t>0$) 
also belongs to the space $C_{-1}(0,+\infty)$, and thus we have to 
distinguish between multiplication with a scalar $\lambda$ in the vector 
space  $C_{-1}^2/\sim$ and  multiplication with a constant 
function in the field $\mathcal{F}_{-1}$ that is defined as follows:
$$
\{\lambda\}\cdot \frac{f}{g} = \frac{\{\lambda\}\, *\,  f}{g},\ \frac{f}{g}\in \mathcal{F}_{-1}.
$$

The unity element of the field $\mathcal{F}_{-1}$ with respect to multiplication is  $I = \frac{\kappa}{\kappa}$ (this equivalence class contains all pairs in the form $(f,\, f)$, $f\in C_{-1}(0,+\infty),\ f\not \equiv 0$).  The unity element $I$ does not belong to the ring $\mathcal{R}_{-1}$ and thus it is not a conventional function.  Indeed, let us suppose that $I$ can be associated with a function $f \in \mathcal{R}_{-1}, \ f\not \equiv 0$. Because of the embedding \eqref{emb}, we then get the following chain of implications:
$$
I = f \ \Leftrightarrow \ \frac{\kappa}{\kappa} = \frac{\kappa * f}{\kappa} \ \Leftrightarrow \  (\kappa * \kappa)(t) = (\kappa * \kappa * f)(t) \ \Rightarrow  \ 
$$
$$
(k^2 * \kappa * \kappa)(t) = (k^2 * \kappa * \kappa * f)(t) \ \Leftrightarrow \ \{1\}^2(t) = (\{1\}^2 * f)(t) \ \Leftrightarrow \ 
t = (t*f)(t) \ \Leftrightarrow \ 
$$
$$
t = (I_{0+}^2\, f)(t) \ \Rightarrow  \ \frac{d^2}{dt^2} t = \frac{d^2}{dt^2} (I_{0+}^2\, f)(t) \ \Leftrightarrow \ 
0 = f(t).
$$
Thus we arrived to a contradiction that means that the unity element $I$ of the field $\mathcal{F}_{-1}$ cannot be reduced to an element of the ring $\mathcal{R}_{-1}$. It can be interpreted as a kind of a generalized function (hyperfunction) that plays the role of the Dirac $\delta$-function in our operational calculus. In fact, most of the convolution quotients from  $\mathcal{F}_{-1}$ are generalized functions (hyperfunctions), not conventional ones. Another important element of this kind is introduced in the following definition:

\begin{definition}
\label{d4}
The element 
\begin{equation}
\label{alg}
S_\kappa = \frac{\kappa}{\kappa^2} \in \mathcal{F}_{-1}
\end{equation}
is called algebraic inverse to the GFI \eqref{GFI}.
\end{definition}
In what follows we denote the field element inverse to an element $e\in \mathcal{F}_{-1}$ with respect to multiplication by $\frac{I}{e}$. 
The background of Definition \ref{d4} is that  on the ring $\mathcal{R}_{-1}$, the GFI \eqref{GFI} can be interpreted as multiplication with the element $\kappa$. The function $\kappa \in \mathcal{R}_{-1}$ corresponds to the element $\frac{\kappa^2}{\kappa}$ of the field $\mathcal{F}_{-1}$.   Thus, the element $S_{\kappa}$ is inverse to $\kappa$ in the field $\mathcal{F}_{-1}$:
$$
{\kappa}\cdot S_\kappa  = \frac{\kappa^2}{\kappa}\cdot \frac{\kappa}{\kappa^2}
= \frac{\kappa^3}{\kappa^3} = \frac{\kappa}{\kappa} = I\ \Leftrightarrow \ S_\kappa = \frac{I}{\kappa}\ \mbox{and} \ \kappa = \frac{I}{S_\kappa}.
$$
The element $S_\kappa$ is a very important object in the framework of the Mikusi\'nski type operational calculus for the GFD with the Sonine kernel $k\in \mathcal{S}_{-1}$ associate to the kernel $\kappa$. In particular, it is employed to define an algebraic general fractional derivative (algebraic GFD): 

\begin{definition}
\label{d5}
Let $k\in \mathcal{S}_{-1}$ be the Sonine kernel associate to the kernel $\kappa$ and  $f\in C[0,+\infty)$. The algebraic GFD of Caputo type  is defined as follows: 
\begin{equation}
\label{AGFD}
 _*\D_{(k)}\, f = S_\kappa \cdot f - S_\kappa \cdot \{ f(0) \},
\end{equation}
where the function $f$ at the right-hand side of the formula  is interpreted as an element of the convolution quotients field $\mathcal{F}_{-1}$.
\end{definition}

It is worth mentioning that the algebraic GFD \eqref{AGFD} is a kind of a generalized derivative that assigns a certain element from $\mathcal{F}_{-1}$ to any function continuous on $[0,+\infty)$, i.e., to any function from the Mikusi\'nski ring $\mathcal{R}$ ($\mathcal{R}\subset \mathcal{R}_{-1}$).  Thus, the algebraic GFD \eqref{AGFD} makes sense also for the functions whose GFD \eqref{GFDC} does not exist in the usual sense. However, for any function  $f\in C_{-1}^1(0,+\infty)$, its GFD exists and coincides with the algebraic GFD \eqref{AGFD}. That's why we use the same notations for both derivatives. 

\begin{theorem}
\label{t9}
For a function $f\in C_{-1}^1(0,+\infty)$, its GFD \eqref{GFDC} exists in the usual sense and coincides with the algebraic GFD \eqref{AGFD}:
\begin{equation}
\label{AGFD-f}
( _*\D_{(k)}\, f)(t) = \mbox{ }_*\D_{(k)}\, f  = S_\kappa \cdot f - S_\kappa \cdot \{ f(0) \}.
\end{equation}
\end{theorem}

\begin{proof}
The inclusion $f\in C_{-1}^1(0,+\infty)$ leads to $f^\prime \in C_{-1}(0,+\infty)$ and $f\in C[0,+\infty)$ (Theorem \ref{t2}). Thus, both the GFD \eqref{GFDC} and the algebraic GFD \eqref{AGFD} exist.  Now we employ the definitions of the derivatives and of the element $S_\kappa$ and then use Theorem \ref{t4} to get the following chain of implications: 
$$
( _*\D_{(k)}\, f)(t) = (\I_{(k)}\, f^\prime)(t) = \phi(t)\in C_{-1}(0,+\infty)\ \Rightarrow 
$$
$$
(\I_{\kappa}\, _*\D_{(k)}\, f)(t) = (\I_{\kappa}\, \phi) (t) \ \Leftrightarrow  \ f(t) - f(0) = (\kappa *\phi)(t) \ \Leftrightarrow  
$$
$$
f(t) = (\kappa *\phi)(t) + f(0) \ \Rightarrow  \ S_\kappa \cdot f  = S_\kappa \cdot (\kappa *\phi) + S_\kappa  \cdot \{ f(0) \} \ \Leftrightarrow 
$$
$$
S_\kappa \cdot f =\frac{\kappa}{\kappa^2}\cdot \frac{ \kappa^2 *\phi}{\kappa} + S_\kappa  \cdot \{ f(0) \} \ \Leftrightarrow \ \frac{\kappa*\phi}{\kappa} = S_\kappa \cdot f - S_\kappa  \cdot \{ f(0) \} \ \Leftrightarrow  
$$
$$
\phi(t) =S_\kappa \cdot f - S_\kappa  \cdot \{ f(0) \}  \ \Leftrightarrow \  ( _*\D_{(k)}\, f)(t) = S_\kappa \cdot f - S_\kappa  \cdot \{ f(0) \}.
$$
\end{proof}

Definition \ref{d5} can be extended to the case of the $n$-fold GFD \eqref{GFDCn}:

\begin{definition}
\label{d6}
Let $k\in \mathcal{S}_{-1}$ be the Sonine kernel associate to the kernel $\kappa$ and  $f\in C^{n-1}[0,+\infty)$. The $n$-fold algebraic GFD of Caputo type is defined as follows: 
\begin{equation}
\label{AGFDn}
 _*\D_{(k)}^n\, f = S_\kappa^n \cdot f - S_\kappa^n \cdot f_n,\ f_n(t) = \sum_{j=0}^{n-1} \, f^{(j)}(0)\, h_{j+1}(t),\ t>0,
\end{equation}
where  the functions $f$ and $f_n$ at the right-hand side of the formula are interpreted as  elements of  $\mathcal{F}_{-1}$.
\end{definition}

Contrary to the $n$-fold  algebraic GFD, the $n$-fold GFD \eqref{GFDCn} exists not for all functions from the space $C^{n-1}[0,+\infty)$. Thus, the $n$-fold  algebraic GFD can be interpreted as a kind of a generalized derivative. However, for the functions 
$f\in C_{-1}^n(0,+\infty)$, Theorem \ref{t7} and the same arguments as the ones we employed in the case $n=1$ lead to validity of the following result:

\begin{theorem}
\label{t10}
For a function $f\in C_{-1}^n(0,+\infty)$, its $n$-fold GFD \eqref{GFDCn} exists in the usual sense and coincides with the $n$-fold algebraic GFD \eqref{AGFDn}:
\begin{equation}
\label{AGFDn-f}
( _*\D_{(k)}^n\, f)(t) = \mbox{ }_*\D_{(k)}^n\, f  = S_\kappa^n \cdot f - S_\kappa^n \cdot f_n,\ f_n(t) = \sum_{j=0}^{n-1} \, f^{(j)}(0)\, h_{j+1}(t),\ t>0.
\end{equation}
\end{theorem}

In the formulas \eqref{AGFD-f} and \eqref{AGFDn-f}, the GFD  \eqref{GFDC}  and the $n$-fold GFD \eqref{GFDCn} are reduced to a simple multiplication in the field $\mathcal{F}_{-1}$ of convolution quotients. These representations allow to transform the initial-value problems for the fractional differential equations containing these general fractional derivatives to some algebraic (in fact, linear) equations in the field $\mathcal{F}_{-1}$ of convolution quotients. In the next section, this technique will be demonstrated  on some examples. However, in general, the solutions to these algebraic equations are the elements of the field $\mathcal{F}_{-1}$, i.e., the generalized functions (hyperfunctions). Because of the embedding \eqref{emb}, in some cases, the elements of $\mathcal{F}_{-1}$ can be interpreted as the conventional functions from the ring  $\mathcal{R}_{-1}$. In the rest of this section, we introduce some important classes of such elements of the field $\mathcal{F}_{-1}$. Our starting point is the following important result:

\begin{theorem} 
\label{t11}
Let $\kappa \in \mathcal{S}_{-1}$ be a Sonine kernel and the convergence radius
$r$ of the power  series
\begin{equation}
\label{ser}
\Sigma(t) = \sum^{+\infty }_{j=0}a_{j}\, t^j,\ a_{j}\in \C,\ t\in \C 
\end{equation}
be non-zero. Then the convolution series 
\begin{equation}
\label{conser}
\Sigma_\kappa(t) = \sum^{+\infty }_{j=0}a_{j}\, \kappa^j(t)
\end{equation}
is convergent for all $t>0$ and defines an  element  of
the ring $\mathcal{R}_{-1}$.
\end{theorem}

\begin{proof}
Because of the inclusion $\kappa \in \mathcal{S}_{-1}$, the kernel $\kappa$ is from the space $C_{-1,0}(0,+\infty)$ and thus it can represented in the form 
\begin{equation}
\label{rep}
\kappa(t) = h_{p}(t)\kappa_1(t),\ t>0,\ 0<p<1,\ \kappa_1 \in C[0,+\infty).
\end{equation}
Let us consider an arbitrary but fixed interval $[0,\, T]$ with $T > 0$. Because the function $\kappa_1(t) = \Gamma(p)t^{1-p}\kappa(t)$ from the representation \eqref{rep} is continuous on $[0,+\infty)$, we get the estimate 
\begin{equation}
\label{est_1}
\exists M_T >0:\ |\kappa_1(t)| = |\Gamma(p)t^{1-p}\kappa(t)| \le M_T,\ t \in [0,\, T].
\end{equation}
Let us fix a point $t_0$, such that $0<|t_0|<r$. Then the series \eqref{ser} is absolutely convergent at the point $t_0$ and we have the following inequalities:
\begin{equation}
\label{est_2}
\exists M >0:\ |a_{j}\, t_0^j| \le M\ \forall j\in \N_0 \ \Rightarrow \ 
|a_{j}| \le \frac{M}{|t_0|^j}\ \forall j\in \N_0.
\end{equation}
To work with the convolution series \eqref{conser}, we need some estimates for the convolution powers $\kappa^j,\ j\ge 2$ on the interval $(0,T]$. Let us start with the case $j=2$:
$$
|\kappa^2(t)|= |(\kappa*\kappa)(t)|\le \int_0^t h_{p}(t-\tau)|\kappa_1(t-\tau)|h_p(\tau)|\kappa_1(\tau)|\, d\tau \le 
$$
$$
M_T^2  (h_{p}*h_{p})(t) = M_T^2\, h_{2p}(t),\ 0< t\le T.
$$
Following the same way, we get the inequalities
\begin{equation}
\label{est_3}
|\kappa^j(t)|\le  M_T^j\, h_{jp}(t) = M_T^j\, \frac{t^{jp-1}}{\Gamma(jp)},\ 0< t\le T,\ j\in \N.
\end{equation}
The inequalities \eqref{est_2} and \eqref{est_3} are now employed to derive the following estimate:
\begin{equation}
\label{est_4}
t^{1-p}|a_j\kappa^j(t)|\le \frac{M}{|t_0|^j} M_T^j\, \frac{t^{(j-1)p}}{\Gamma(jp)} \le 
\frac{M}{T^p}\frac{\left(\frac{M_T\, T^p}{|t_0|}\right)^j}{\Gamma(jp)},\ j\in \N,\ 0\le t\le T.
\end{equation}
Because of the Stirling asymptotic formula
$$
\Gamma(x+1) \sim \sqrt{2\pi x}\left( \frac{x}{e}\right)^x,\ x\to +\infty,
$$
the number series 
$$
\sum_{j=0}^{+\infty}\frac{M}{T^p}\frac{\left(\frac{M_T\, T^p}{|t_0|}\right)^j}{\Gamma(jp)}
$$
is absolutely convergent. The estimate \eqref{est_4} allows us to employ the Weierstrass M-test that says that the series
\begin{equation}
\label{conser_1}
t^{1-p}\sum_{j=0}^{+\infty}a_j\kappa^j(t)
\end{equation}
is absolutely and uniformly convergent on the interval $[0,\, T]$. Furthermore, the functions $t^{1-p}a_j\kappa^j(t),\ j\in \N_0$ are continuous on $[0,\, T]$ (see the inequality \eqref{est_3} and remember that $p\in (0,1)$) and thus the uniform limit theorem ensures that the series \eqref{conser_1} is a function continuous on the interval $[0,\, T]$. Because $T$ can be chosen arbitrary large, this means that the convolution series \eqref{conser}
is convergent for all $t>0$ and defines an  element  of
the ring $\mathcal{R}_{-1}$.
\end{proof}

It is worth mentioning that the convolution series \eqref{conser} is a far reaching generalization of the power series \eqref{ser} that corresponds to the case of the kernel $\kappa = \{1\}$ (Mikusi\'nski's operational calculus for the first order derivative).
As an example, we consider the geometric series
\begin{equation}
\label{geom}
\Sigma(t) = \sum_{j=1}^{+\infty} \lambda^{j-1}t^j,\ \lambda \in \C,\ t\in \C.
\end{equation}
If $\lambda=0$, the series \eqref{geom} contains just one term: $\Sigma(t)=t$. For $\lambda \not =0$, its convergence radius $r = 1/|\lambda|$ is non-zero. According to Theorem \eqref{t11}, the convolution series ($\kappa \in \mathcal{S}_{-1}$) 
\begin{equation}
\label{l}
l_{\kappa,\lambda}(t) = \sum_{j=1}^{+\infty} \lambda^{j-1}\kappa^j(t),\ \lambda \in \C
\end{equation}
is convergent for all $t>0$ and defines a function from the space $C_{-1}(0,+\infty)$. 

In the framework of the  Mikusi\'nski operational calculus (\cite{Mik59}), the kernel function $\kappa$ is the constant function $\{ 1\}$ and thus $\kappa^j(t) = \{ 1\}^j(t) = h_{j}(t)$. Then the convolution series \eqref{l} becomes a familiar power series for the exponential function:
\begin{equation}
\label{l-Mic}
l_{\kappa,\lambda}(t) = \sum_{j=1}^{+\infty} \lambda^{j-1}h_j(t) = 
\sum_{j=0}^{+\infty} \frac{(\lambda\, t)^j}{j!} = e^{\lambda\, t}.
\end{equation}

The operational calculi for the Riemann-Liouville fractional derivative (\cite{HadLuc}, \cite{LucSri95}) and for  the Caputo fractional derivative (\cite{LucGor99}) are based on  the kernel $\kappa = h_{\alpha}$ of the Riemann-Liouville fractional integral. In this case $\kappa^j(t) = h_{\alpha}^j(t) = h_{j\alpha}(t)$ and the convolution series \eqref{l} takes the form
\begin{equation}
\label{l-Cap}
l_{\kappa,\lambda}(t) = \sum_{j=1}^{+\infty} \lambda^{j-1}h_{j\alpha}(t) = 
t^{\alpha-1}\sum_{j=0}^{+\infty} \frac{\lambda^j\, t^{j\alpha}}{\Gamma(j\alpha+\alpha)} = t^{\alpha -1}E_{\alpha,\alpha}(\lambda\, t^{\alpha}),
\end{equation}
where the two-parameters Mittag-Leffler function $E_{\alpha,\alpha}$ is defined by \eqref{ML}. 

Another interesting case is the GFD \eqref{GFD_3} that corresponds to the GFI \eqref{GFI_3} with the kernel $\kappa(t) = h_{\alpha,\rho}(t)=h_\alpha(t)\, e^{-\rho t},\ \ \alpha >0,\ \rho \ge 0,\ t>0$. According to the formula \eqref{2-10}, $\kappa^j(t) = 
h_{j\alpha,\rho}(t)$ and we get the following representation of the convolution series \eqref{l}:
\begin{equation}
\label{l-exp}
l_{\kappa,\lambda}(t) = \sum_{j=1}^{+\infty} \lambda^{j-1}h_{j\alpha,\rho}(t) = 
e^{-\rho t}t^{\alpha -1}\sum_{j=0}^{+\infty} \frac{\lambda^j\, t^{j\alpha}}{\Gamma(j\alpha+\alpha)} = e^{-\rho t}t^{\alpha -1}E_{\alpha,\alpha}(\lambda\, t^{\alpha}).
\end{equation}

In the framework of the operational calculus for the GFD \eqref{GFD_2} with the  Mittag-Leffler function \eqref{3-9} in the kernel, the corresponding GFI \eqref{GFI_2} has the kernel $\kappa(t) = h_{1-\beta+\alpha}(t)\, +\, h_{1-\beta}(t),\ 0<\alpha < \beta <1$. The convolution powers of this kernel are given by the formula \eqref{2-16} and the convolution series \eqref{l} takes the following form:
\begin{equation}
\label{l-ML-d}
l_{\kappa,\lambda}(t) = \sum_{j=1}^{+\infty} \lambda^{j-1} \sum_{i=0}^j \binom{j}{i} h_{i\alpha + j(1-\beta)}(t) = 
\end{equation}
$$
\frac{1}{\lambda t} \sum_{j=0}^{+\infty}\sum_{l_1+l_2 =j} \frac{j!}{l_1!l_2!}\frac{(\lambda t^{1-\beta+\alpha})^{l_1}(\lambda t^{1-\beta})^{l_2}}{\Gamma(l_1(1-\beta+\alpha)+l_2(1-\beta))} = 
$$
$$
\frac{1}{\lambda t}  E_{(1-\beta,1-\beta+\alpha),0}(\lambda t^{1-\beta},\lambda t^{1-\beta+\alpha}),
$$
where $E_{(1-\beta,1-\beta+\alpha),0}$ is a particular case of the multinomial Mittag-Leffler function
\begin{equation}
\label{MLm}
E_{(\alpha_1,\ldots,\alpha_m),\beta}(z_1,\ldots,z_m):=
\sum_{j=0}^{+\infty} \sum_{l_1+\cdots +l_m =j}\frac{j!}{l_1!\times\cdot\times
l_m!}\frac{\prod_{i=1}^m  
z_i^{l_i}}{\Gamma(\beta+\sum_{i=1}^m \alpha_i l_i)}
\end{equation}
introduced for the first time in \cite{Luc93} (see also \cite{HadLuc}).

The function $l_{\kappa,\lambda}$ defined by the convolution series \eqref{l} plays a very important role in applications of the constructed operational calculus for analytical treatment of the fractional differential equations with the GFDs. The meaning of $l_{\kappa,\lambda}$ is explained in the following theorem:

\begin{theorem}
\label{t12}
For any $\lambda \in \C$, the element $l_{\kappa,\lambda} \in \mathcal{R}_{-1}\subset \mathcal{F}_{-1}$  defined by the convolution series \eqref{l} is  inverse to the element $(S_\kappa -\lambda) \in \mathcal{F}_{-1}$, i.e., the following operational relation hold true:
\begin{equation}
\label{op-rel}
\frac{I}{S_{\kappa} - \lambda}\, = \, l_{\kappa,\lambda}(t),\ t>0.
\end{equation}
\end{theorem}

\begin{proof} 
If $\lambda = 0$, the convolution series \eqref{l} reduces to a single term: $l_{\kappa,\lambda}= \kappa$. The relation \eqref{op-rel} takes then the form $\frac{I}{S_{\kappa}} = \kappa$ that holds true because of Definition \ref{d4}. Now we consider the case $\lambda \not = 0$ and start with the evident formula
\begin{equation}
\label{rel}
\kappa\cdot(S_\kappa -\lambda) = \kappa\cdot S_\kappa - \lambda \kappa = I - \lambda\kappa.
\end{equation}
Then we proceed with calculation of the following product:
$$
(I - \lambda\kappa)\cdot l_{\kappa,\lambda} = l_{\kappa,\lambda} - \lambda (\I_{(\kappa)}\, \sum_{j=1}^{+\infty} \lambda^{j-1}\kappa^j)(t) = l_{\kappa,\lambda} - \lambda  \sum_{j=1}^{+\infty} \lambda^{j-1} (\I_{(\kappa)}\, \kappa^j)(t) =
$$
$$
l_{\kappa,\lambda} - \sum_{j=1}^{+\infty} \lambda^{j}  \kappa^{j+1}(t) =
\sum_{j=1}^{+\infty} \lambda^{j-1}\kappa^j(t) - \sum_{j=2}^{+\infty} \lambda^{j-1}\kappa^j(t) = \kappa.
$$
In the formula above, we could apply the fractional integral $\I_{(\kappa)}$ to the convolution series defining the element $l_{\kappa,\lambda}$ term by term because this series is absolutely and uniformly convergent on any interval $[\epsilon,T],\ 0<\epsilon <T$ (see the proof of Theorem \eqref{t11}). 

Combining the last formula with the relation \eqref{rel}, we arrive at the representation
$$
\kappa\cdot(S_\kappa -\lambda)\cdot l_{\kappa,\lambda} = (I - \lambda\kappa)\cdot l_{\kappa,\lambda} = \kappa,
$$
which implicates the identity
$$
(S_\kappa -\lambda)\cdot l_{\kappa,\lambda} = I
$$
that proves the operational relation \eqref{op-rel}. 
\end{proof}

Making use of the formulas \eqref{l-Mic}-\eqref{l-ML-d}, we list now the following important operational relations, both the known and the new ones:

1) $\kappa(t) \equiv  1,\ t\ge0$ (Mikusi\'nski's operational calculus \cite{Mik59} for the first order derivative):
\begin{equation}
\label{l-Mic-op}
\frac{I}{S_{\kappa} - \lambda}\, = \,   e^{\lambda\, t}.
\end{equation}

2) $\kappa(t) = h_{\alpha}(t),\ t>0$ (operational calculi for the Riemann-Liouville and  the Caputo fractional derivatives \cite{HadLuc,Luc93,LucGor99,LucSri95}):
\begin{equation}
\label{l-Cap-op}
\frac{I}{S_{\kappa} - \lambda}\, = \, t^{\alpha -1}E_{\alpha,\alpha}(\lambda\, t^{\alpha}),
\end{equation}
where the two-parameters Mittag-Leffler function $E_{\alpha,\alpha}$ is defined by \eqref{ML}. 

3)  $\kappa(t) = h_{\alpha,\rho}(t)=h_\alpha(t)\, e^{-\rho t},\ \ \alpha >0,\ \rho\ge 0,\ t>0$ (operational calculus for the GFD in form \eqref{GFD_3}):  
\begin{equation}
\label{l-exp-op}
\frac{I}{S_{\kappa} - \lambda}\, = \, e^{-\rho t}t^{\alpha -1}E_{\alpha,\alpha}(\lambda\, t^{\alpha}).
\end{equation}

4) $\kappa(t) = h_{1-\beta+\alpha}(t)\, +\, h_{1-\beta}(t),\ 0<\alpha < \beta <1$  (operational calculus for the GFD \eqref{GFD_2} with the  Mittag-Leffler function \eqref{3-9} in the kernel): 
\begin{equation}
\label{l-ML-d-op}
\frac{I}{S_{\kappa} - \lambda}\, = \, \frac{1}{\lambda t}  E_{(1-\beta,1-\beta+\alpha),0}(\lambda t^{1-\beta},\lambda t^{1-\beta+\alpha}),
\end{equation}
where the function $E_{(1-\beta,1-\beta+\alpha),0}$ is a particular case of the multinomial Mittag-Leffler function \eqref{MLm}.

Because the right-hand side of the operational formula \eqref{op-rel} is a function from the ring $\mathcal{R}_{-1}$, we can deduce other useful operational relations by calculating the convolution powers of the function $l_{\kappa,\lambda}$. The general formula takes then the form
\begin{equation}
\label{op-rel-m}
\frac{I}{(S_{\kappa} - \lambda)^m}\, = \, l^m_{\kappa,\lambda}(t),\ t>0,\ m\in \N.
\end{equation}
The convolution powers $l^m_{\kappa,\lambda}$ can be calculated in explicit form. For instance, for $m=2$, we can build the Cauchy product for the convolution series defining the element $l_{\kappa,\lambda}$ (this series is absolutely and uniformly convergent on any interval $[\epsilon,T],\ 0<\epsilon<T$) and thus obtain the formula
\begin{equation}
\label{op-rel-2}
\frac{I}{(S_{\kappa} - \lambda)^2}\, = \, 
\sum_{j=2}^{+\infty} (j-1)\lambda^{j-2}\kappa^j(t)
\, = \, 
\left( \kappa\, *\, \sum_{j=2}^{+\infty} (j-1)\lambda^{j-2}\kappa^{j-1}\right)(t).
\end{equation}
The same procedure can be applied to derive a more general operational relation
\begin{equation}
\label{Cauchy}
\frac{I}{(S_{\kappa} - \lambda)^m}\, = \,  
\sum_{j=m}^{+\infty} a_{j,m}\lambda^{j-m}\kappa^j(t)
\, = \, 
\left( \kappa^{m-1}\, *\, \sum_{j=m}^{+\infty} a_{j,m}\lambda^{j-m}\kappa^{j-m+1}\right)(t),
\end{equation}
where the coefficients $a_{j,m}$ can be calculated in explicit form (see the formula \eqref{op-rel-2}). It is worth mentioning that the series from the right-hand side of the formula \eqref{Cauchy} belongs to the space $C_{-1}(0,+\infty)$ according to Theorem \ref{t11}. 

For some particular cases of the kernel function $\kappa$, the convolution powers $l^m_{\kappa,\lambda}$ can be calculated as $m$-fold convolutions of $l_{\kappa,\lambda}$. For instance, in the case $\kappa(t)\equiv  1,\ t\ge0$ (Mikusi\'nski's operational calculus for the first order derivative), we easily get the well-known operational relation (\cite{Mik59}):
\begin{equation}
\label{l-Mic-op-m}
\frac{I}{(S_{\kappa} - \lambda)^m}\, = \,  h_m(t)\, e^{\lambda\, t}.
\end{equation} 
For $\kappa(t) = h_{\alpha}(t),\ t>0$ (operational calculus for  the Caputo fractional derivative), the operational relation \eqref{op-rel-m} takes the form (\cite{LucGor99}):
\begin{equation}
\label{l-Cap-op-m}
\frac{I}{(S_{\kappa} - \lambda)^m}\, = \, t^{m\alpha -1}E_{\alpha,m\alpha}^m(\lambda t^\alpha),\ t>0,\ m\in \N, 
\end{equation}
where the Mittag-Leffler type function $E_{\alpha,\beta}^{m}$ is defined as follows:
$$
E_{\alpha,\beta}^{m}(z):=
\sum^{\infty
}_{j=0}\frac{(m)_j z^{j}}{  j!
\Gamma (\alpha j + \beta)}, \ \alpha,\beta>0, \ z\in \C, \ 
(m)_j = \prod_{i=0}^{j-1} (m+i).
$$

\vspace*{-3pt} 

\section{Fractional differential equations with the GFDs}
\label{sec:6}

\setcounter{section}{5}
\setcounter{equation}{0}\setcounter{theorem}{0}

In this section, we demonstrate the operational method for analytical treatment of the initial-value problems for the fractional differential equations with the GFDs with the Sonine kernels from $\mathcal{S}_{-1}$ on some examples. 

Let us start with the following simple but still not completely investigated case:
\begin{equation}
\label{eq-1-1}
\begin{cases}
( _*\D_{(k)}\, y)(t) - \lambda y(t) = f(t), & \lambda \in \R,\ t>0, \\
y(0) = y_0, & y_0\in \R .
\end{cases}
\end{equation}
In the equation \eqref{eq-1-1}, $_*\D_{(k)}$ is the Caputo type GFD with a kernel $k \in \mathcal{S}_{-1}$, the source function $f$ is from the space $C_{-1}^1(0,+\infty)$ and the unknown function $y$ is looked for in the space $C_{-1}^1(0,+\infty)$. 

Because of the inclusion $y\in C_{-1}^1(0,+\infty)$, we can apply the formula \eqref{AGFD-f} from Theorem \ref{t9} and transform the initial-value problem \eqref{eq-1-1} to a single linear equation in the field $\mathcal{F}_{-1}$ of convolution quotients:
\begin{equation}
\label{eq-1-2}
S_\kappa \cdot y - S_\kappa \cdot \{ y_0 \} - \lambda y = f, 
\end{equation}
where the functions $y,\, f \in \mathcal{R}_{-1}$ are interpreted as elements of $\mathcal{F}_{-1}$. 

The linear equation \eqref{eq-1-2} always has a unique solution in the field $\mathcal{F}_{-1}$  that can be represented in the following form:
\begin{equation}
\label{eq-1-3}
y  = f \cdot \frac{I}{S_\kappa - \lambda} +  (S_\kappa \cdot \{ y_0 \}) \cdot \frac{I}{S_\kappa - \lambda}.
\end{equation}

Thus, the formula \eqref{eq-1-3} provides us with the unique solution to the initial-value problem \eqref{eq-1-1}. However, it is an element from the field $\mathcal{F}_{-1}$ of convolution quotients, i.e., a generalized function (hyperfunction). Luckily, this generalized function can be interpreted as an element of the ring $\mathcal{R}_{-1}$, i.e., as a conventional function:
\begin{equation}
\label{eq-1-4}
y(t)  = y_f(t) + y_{iv}(t),\ y_f:= f \cdot \frac{I}{S_\kappa - \lambda},\ y_{iv}:= (S_\kappa \cdot \{ y_0 \}) \cdot \frac{I}{S_\kappa - \lambda}.
\end{equation}

Indeed, the operational relation \eqref{op-rel} from Theorem \ref{t12} represents the element $\frac{I}{S_{\kappa} - \lambda}$ of the field $\mathcal{F}_{-1}$ as the conventional function $l_{\kappa,\lambda} \in \mathcal{R}_{-1}$ in form of the convolution series \eqref{l}. Embedding of the ring $\mathcal{R}_{-1}$ into the field $\mathcal{F}_{-1}$ immediately leads to the representation
\begin{equation}
\label{eq-1-5}
y_f (t) = f \cdot \frac{I}{S_\kappa - \lambda} = (f\, *\, l_{\kappa,\lambda})(t) = 
\int_0^t l_{\kappa,\lambda}(\tau)f(t-\tau)\, d\tau.
\end{equation}
As to the part $y_{iv}$ of the solution \eqref{eq-1-4}, we first deduce the following important operational relation:
\begin{equation}
\label{op-rel-1}
(S_\kappa \cdot \{ 1 \}) \cdot \frac{I}{S_\kappa - \lambda} \, = \, L_{\kappa,\lambda}(t),\ t>0,
\end{equation}
where the function $L_{\kappa,\lambda} \in \mathcal{R}_{-1}$ is defined in form of a convolution series:
\begin{equation}
\label{L}
L_{\kappa,\lambda}(t) = (k \, *\ l_{\kappa,\lambda})(t) = 1 + \{ 1 \} * \sum_{j=1}^{+\infty} \lambda^{j}\kappa^j(t),\ \lambda \in \C.
\end{equation}
Indeed, first we mention the formula
\begin{equation}
\label{op-rel-3}
S_\kappa \cdot \{ 1 \} = \frac{\kappa}{\kappa*\kappa}\cdot \frac{ \{1\}*\kappa}{\kappa} =
\frac{\{1\}*\kappa^2}{\kappa^3} = \frac{\{1\}}{\kappa} = k,
\end{equation}
where the kernel $k\in \mathcal{S}_{-1}$ is the Sonine kernel associated to the kernel $\kappa$ (see the Sonine condition \eqref{3-2}). Then we again employ Theorem \ref{t12} and the operational relation \eqref{op-rel} to deduce the inclusion $L_{\kappa,\lambda} \in \mathcal{R}_{-1}$ (as a convolution of two functions from $\mathcal{R}_{-1}$) and  the representation \eqref{L}:
$$
(S_\kappa \cdot \{ 1 \}) \cdot \frac{I}{S_\kappa - \lambda} = (k \, *\ l_{\kappa,\lambda})(t) = (k * \sum_{j=1}^{+\infty} \lambda^{j-1}\kappa^j)(t) =
$$
$$
 (\I_{(k)}\, \sum_{j=1}^{+\infty} \lambda^{j-1}\kappa^j)(t) = \sum_{j=1}^{+\infty} \lambda^{j-1}(\I_{(k)}\, \kappa^j)(t) 
=
 1 + \{ 1 \} * \sum_{j=1}^{+\infty} \lambda^{j}\kappa^j(t).
$$
In the above derivations, we could apply the fractional integral $\I_{(k)}$ to the convolution series defining the element $l_{\kappa,\lambda}$ term by term because this series is absolutely and uniformly convergent on any interval $[\epsilon,T],\ 0<\epsilon <T$ (see the proof of Theorem \eqref{t11}). 

The  results presented above are summarized in the following theorem:
\begin{theorem}
\label{t13}
Let $k\in \mathcal{S}_{-1}$ be a Sonine kernel.  The unique solution $y\in C_{-1}^1(0,+\infty)$  to the initial-value problem \eqref{eq-1-1} for the fractional differential equation with the GFD of Caputo type can be represented in the form
\begin{equation}
\label{sol-1}
y (t) = (f\, *\, l_{\kappa,\lambda})(t) + y_0\,  L_{\kappa,\lambda}(t),
\end{equation}
where $\kappa\in \mathcal{S}_{-1}$ is the Sonine kernel  associate 
to  the kernel $k$ and $l_{\kappa,\lambda}$, $L_{\kappa,\lambda}$ are the convolution series \eqref{l} and \eqref{L}, respectively. 
\end{theorem}

To finalize the proof of this theorem, we just mention that the inclusion $y\in C_{-1}^1(0,+\infty)$ of the solution \eqref{sol-1} follows from Theorem \ref{t1} and the representations \eqref{l} and \eqref{L} of the functions $l_{\kappa,\lambda}$ and $L_{\kappa,\lambda}$. Furthermore, the part $y=y_f(t) = f\, *\, l_{\kappa,\lambda}(t)$ of the solution formula \eqref{sol-1} evidently solves the initial-value problem \eqref{eq-1-1} for the inhomogeneous equation and homogeneous initial condition ($y_0 = 0$) whereas the function $y = y_{iv}(t) = y_0\,  L_{\kappa,\lambda}(t)$ solves the problem \eqref{eq-1-1} for the homogeneous equation ($f(t) \equiv 0,\ t>0$) and the inhomogeneous initial condition. 

Now we discuss some particular cases of the initial-value problem \eqref{eq-1-1}, both the known and the new ones. 

1) $\kappa(t) \equiv  1,\ t\ge0$ (Mikusi\'nski's operational calculus  for the first order derivative). The function $l_{\kappa,\lambda}$ is given by the formula \eqref{l-Mic-op} and 
\begin{equation}
\label{L-Mic-op}
L_{\kappa,\lambda}(t) = 1 + \{ 1 \} * \sum_{j=1}^{+\infty} \lambda^{j}h_{j}(t) = 
1 + \sum_{j=1}^{+\infty} \lambda^{j}h_{j+1}(t) = e^{\lambda t}.
\end{equation}
According to Theorem \ref{t13}, the initial-value problem 
\begin{equation}
\label{eq-1-1-1}
\begin{cases}
y^\prime(t) - \lambda y(t) = f(t), & t>0, \\
y(0) = y_0, & y_0\in \R .
\end{cases}
\end{equation}
has the well-known solution 
\begin{equation}
\label{sol-1-1}
y (t) = (f\, *\, e^{\lambda \tau})(t) + y_0\, e^{\lambda\, t}.
\end{equation}

2) $\kappa(t) = h_{\alpha}(t),\ t>0,\ 0<\alpha <1$ (operational calculus for  the Caputo fractional derivative \eqref{CD}). The function $l_{\kappa,\lambda}$ is given by the formula \eqref{l-Cap-op} and the function $L_{\kappa,\lambda}$ can be represented as follows:
\begin{equation}
\label{L-Cap-op}
L_{\kappa,\lambda}(t)\, = \, 1 + \{ 1 \} * \sum_{j=1}^{+\infty} \lambda^{j}h_{j\alpha}(t) = 
1 + \sum_{j=1}^{+\infty} \lambda^{j}h_{j\alpha+1}(t) = E_{\alpha,1}(\lambda\, t^{\alpha}),
\end{equation}
where the Mittag-Leffler function $E_{\alpha,1}$ is defined by \eqref{ML}. 
The initial-value problem
\begin{equation}
\label{eq-1-1-2}
\begin{cases}
( _*D^\alpha\, y)(t) - \lambda y(t) = f(t), & t>0, \ 0<\alpha < 1,\\
y(0) = y_0, & y_0\in \R 
\end{cases}
\end{equation}
for the fractional differential equation with the Caputo derivative $ _*D^\alpha$  has then the well-known solution (\cite{LucGor99})
\begin{equation}
\label{sol-1-2}
y (t) = (f\, *\, \tau^{\alpha -1}E_{\alpha,\alpha}(\lambda\, \tau^{\alpha}))(t) + y_0\, E_{\alpha,1}(\lambda\, t^{\alpha}).
\end{equation}

3)  $\kappa(t) = h_{\alpha,\rho}(t)=h_\alpha(t)\, e^{-\rho t},\ \ \alpha >0,\ \rho \ge 0,\ t>0$ (operational calculus for the GFD in form \eqref{GFD_3}). In this case, the function $l_{\kappa,\lambda}$ is given by the formula \eqref{l-exp-op} and the function $L_{\kappa,\lambda}$ is as follows:
$$
L_{\kappa,\lambda}(t)\, = \, (k\, * l_{\kappa,\lambda})(t) = \left( (h_{1-\alpha,\rho} + \rho\{1\}*h_{1-\alpha,\rho}), \sum_{j=1}^{+\infty} \lambda^{j-1}h_{j\alpha,\rho}\right)(t) = 
$$
$$
\sum_{j=1}^{+\infty} \lambda^{j-1}h_{(j-1)\alpha+1,\rho}(t) + \rho\left(\{1\}* 
\sum_{j=1}^{+\infty} \lambda^{j-1}h_{(j-1)\alpha+1,\rho}\right)(t) = 
$$
\begin{equation}
\label{L-exp-op}
e^{-\rho t}E_{\alpha,1}(\lambda\, t^{\alpha}) + \rho \int_0^t e^{-\rho \tau}E_{\alpha,1}(\lambda\, \tau^{\alpha})\, d\tau.
\end{equation} 
Specializing Theorem \ref{t13} to this case, we get the following result:

The initial-value problem \eqref{eq-1-1} with the GFD  
$$
( _*\D_{(k)}\, f)(t) = \int_0^t \left( h_{1-\alpha,\rho} (t-\tau) \, +\, \frac{\rho^\alpha}{\Gamma(1-\alpha)}\, \gamma(1-\alpha,\rho (t-\tau)) \right)\, f^\prime(\tau)\, d\tau
$$
has a unique solution given by the formula
\begin{equation}
\label{sol-1-3}
y (t) = (f\, *\, e^{-\rho \tau}\tau^{\alpha -1}E_{\alpha,\alpha}(\lambda\, \tau^{\alpha})(t) + 
\end{equation}
$$
y_0\, \left(e^{-\rho t}E_{\alpha,1}(\lambda\, t^{\alpha}) + \rho \int_0^t e^{-\rho \tau}E_{\alpha,1}(\lambda\, \tau^{\alpha})\, d\tau\right).
$$

4) Finally we consider the case of the kernel $\kappa(t) = h_{1-\beta+\alpha}(t)\, +\, h_{1-\beta}(t),\ 0<\alpha < \beta <1$  (operational calculus for the GFD \eqref{GFD_2} with the  Mittag-Leffler function \eqref{3-9} in the kernel). The function $l_{\kappa,\lambda}$ is given by the formula  \eqref{l-ML-d-op}. To get an explicit expression for the function $L_{\kappa,\lambda}$, we put the representation \eqref{2-16} into the formula \eqref{L} and arrive at the following result:
$$
L_{\kappa,\lambda}(t)\, = \, \, 1 + \{ 1 \} * \sum_{j=1}^{+\infty} \lambda^{j}\kappa^j(t) = 
$$
\begin{equation}
\label{L-ML-d-op}
\, 1 + \{ 1 \} * \sum_{j=1}^{+\infty} \lambda^{j} \sum_{i=0}^j \binom{j}{i} h_{i\alpha + j(1-\beta)}(t) = 
 E_{(1-\beta,1-\beta+\alpha),1}(\lambda t^{1-\beta},\lambda t^{1-\beta+\alpha}),
\end{equation}
where the function $E_{(1-\beta,1-\beta+\alpha),1}$ is a particular case of the multinomial Mittag-Leffler function \eqref{MLm}. 
According to Theorem \ref{t13}, the unique solution to the initial-value problem \eqref{eq-1-1} with the GFD  
$$
( _*\D_{(k)}\, f)(t) = \int_0^t \tau^{\beta -1}\, E_{\alpha,\beta}(-\tau^\alpha)\,   f^\prime(t-\tau)\, d\tau,\ 0<\alpha < \beta <1,\ t>0
$$
is given by the formula
\begin{equation}
\label{sol-1-4}
y (t) = (f\, *\, \frac{1}{\lambda \tau}  E_{(1-\beta,1-\beta+\alpha),0}(\lambda \tau^{1-\beta},\lambda \tau^{1-\beta+\alpha}))(t) +
\end{equation}
$$
y_0\,  E_{(1-\beta,1-\beta+\alpha),1}(\lambda t^{1-\beta},\lambda t^{1-\beta+\alpha}).
$$

Another important observation is that the initial-value problem \eqref{eq-1-1} for the homogeneous fractional differential equation and with the initial condition $y(0)=y_0 =1$ coincides with the problem \eqref{relax} considered in \cite{Koch11}. Thus, we have the following result:

\begin{theorem}
\label{t14}
Let the kernel $k \in \mathcal{S}_{-1}$ of the GFD \eqref{GFDC} belong to the set $\mathcal{K}$ of the Sonine kernels, i.e., the conditions K1)-K4) (see Introduction) hold true, and $\kappa \in \mathcal{S}_{-1}$ be its associate kernel. 

Then, for $\lambda <0$, the convolution series $L_{\kappa,\lambda}$ defined by \eqref{L} is continuous on $[0,\, +\infty)$ and infinitely differentiable and completely monotone on $\R_+$. 
\end{theorem}
Indeed, according to Theorem \ref{t13}, the convolution series $L_{\kappa,\lambda}$ is the unique solution to the initial-value problem \eqref{relax}. For $\lambda <0$ and $k \in \mathcal{K}$ we can employ the result proved in \cite{Koch11} (see the problem (A) in the Introduction) that immediately leads to the statement of Theorem \ref{t14}. 

It is an easy exercise to verify that the kernels of the GFDs considered in the examples 1), 2), and 4) above, satisfy the conditions K1)-K4) and thus we can apply Theorem \ref{t14} to the convolution series \eqref{L-Mic-op}, \eqref{L-Cap-op},  and \eqref{L-ML-d-op}. For $\lambda <0$, the first two functions, $L_{\kappa,\lambda}(t) = e^{\lambda\, t}$ and $L_{\kappa,\lambda}(t) = E_{\alpha,1}(\lambda\, t^{\alpha})$ are known to be completely monotone, see \cite{[SSV]}. However, the statement that the function
$$
L_{\kappa,\lambda}(t) = 
 E_{(1-\beta,1-\beta+\alpha),1}(\lambda t^{1-\beta},\lambda t^{1-\beta+\alpha}),\ \lambda <0,
$$
is completely monotone seems to be a new result ($E_{(1-\beta,1-\beta+\alpha),1}$ is a particular case of the multinomial Mittag-Leffler function \eqref{MLm}). 

Now we consider a generalization of the initial-value problem \eqref{eq-1-1} to the case of a linear multi-term fractional differential equation:
\begin{equation}
\label{eq-2-1}
\begin{cases}
\sum_{j=0}^{n} a_j( _*\D_{(k)}^j\, y)(t) = f(t), & a_j\in \R,\ a_n \not = 0,\, t>0, \\
y^{(j)}(0) = y_{j,0}, \ j=0,\dots,n-1, & y_{j,0}\in \R .
\end{cases}
\end{equation}
In \eqref{eq-2-1}, $_*\D_{(k)}^j$ is the $j$-fold GFD \eqref{GFDCn} with a kernel $k \in \mathcal{S}_{-1}$, the source function $f$ is from the space $C_{-1}^1(0,+\infty)$ and the unknown function $y$ is looked for in the space $C_{-1}^n(0,+\infty)$. 

Similar to the case of the initial-value problem \eqref{eq-1-1},  for $j\ge 1$ and $y\in C_{-1}^n(0,+\infty)$, we employ the formula \eqref{AGFDn-f} from Theorem \ref{t10} and thus transform the initial-value problem \eqref{eq-2-1} to a single linear equation in the field $\mathcal{F}_{-1}$ of convolution quotients:
\begin{equation}
\label{eq-2-2}
a_0\, y + \sum_{j=1}^{n} a_j \left(S_\kappa^j \cdot y - S_\kappa^j \cdot y_j \right) = f, \ y_j(t) = \sum_{i=0}^{j-1} \, y_{i,0}\, h_{i+1}(t),
\end{equation}
where the functions $y,\, y_j,\ f \in \mathcal{R}_{-1}$ are interpreted as the elements of $\mathcal{F}_{-1}$. 

The unique solution to the linear equation \eqref{eq-2-2} in the field $\mathcal{F}_{-1}$ can be represented as follows:
\begin{equation}
\label{eq-2-3}
y  = f \cdot \frac{I}{a_0 + \sum_{j=1}^{n} a_j S_\kappa^j} +  \frac{\sum_{j=1}^{n} a_j S_\kappa^j \cdot y_j } {a_0 + \sum_{j=1}^{n} a_j S_\kappa^j}.
\end{equation}

Now we show that the generalized solution \eqref{eq-2-3} can be interpreted as an element of the ring $\mathcal{R}_{-1}$, i.e., as a conventional function:
\begin{equation}
\label{eq-2-4}
y(t)  = y_f(t) + y_{iv}(t),\ y_f:= f \cdot \frac{I}{P_n(S_\kappa)},\ y_{iv}:= \frac{P_{iv}(S_\kappa)} {P_n(S_\kappa)},
\end{equation}
where the polynomial $P_n(S_\kappa)$ and the function $P_{iv}(S_\kappa)$ are defined as follows:
\begin{equation}
\label{pn}
P_n(S_\kappa) := a_0 + \sum_{j=1}^{n} a_j S_\kappa^j, \ 
P_{iv}(S_\kappa) := \sum_{j=1}^{n} a_j S_\kappa^j \cdot y_j.
\end{equation}
To get an explicit formula for the function $y_f$, we employ a representation of  the 
rational function $I/P_n(S_\kappa)$ in the variable $S_\kappa$  as a sum of the partial fractions on the field $\mathcal{F}_{-1}$:
\begin{equation}
\label{y-f-1}
\frac{I}{P_n(S_\kappa)} = \sum_{i=1}^l \sum_{j=1}^{p_i} c_{ij}\frac{I}{(S_{\kappa}-\lambda_i)^j},\ p_i\in \N,\ \sum_{i=1}^l p_i = n.
\end{equation}
Using this formula and the operational relation \eqref{op-rel-m}, we get the representation
\begin{equation}
\label{y-f-2}
G_{\kappa,n}(t) := \frac{I}{P_n(S_\kappa)} = \sum_{i=1}^l \sum_{j=1}^{p_i} c_{ij}l_{\kappa,\lambda_i}^j(t),
\end{equation}
where the function $l_{\kappa,\lambda}$ is defined by the convolution series \eqref{l}. 
Moreover, 
the formula \eqref{Cauchy} allows us to represent the function $G_{\kappa,n},\ n>1$ as follows:
\begin{equation}
\label{y-f-3}
G_{\kappa,n}(t) = \frac{I}{a_n\prod_{i=1}^l (S_\kappa - \lambda_i)^{p_i}} = (\kappa^{n-1}\, * \, \tilde{G}_{\kappa,n})(t),\ \tilde{G}_{\kappa,n} \in \mathcal{R}_{-1}.
\end{equation}
Thus the function $y_f$ takes the form
\begin{equation}
\label{y-f-4}
y_f(t) = (f\, * \, G_{\kappa,n})(t) = \int_0^t G_{\kappa,n}(\tau)f(t-\tau)\, d\tau,
\end{equation}
where  $G_{\kappa,n}$ is defined by the formula \eqref{y-f-2}.

Now we handle the part $y_{iv}$ of the solution \eqref{eq-2-4} and first employ the formula  $S_\kappa \cdot \{ 1 \}  = k$ (see \eqref{op-rel-3}) to rewrite the numerator $P_{iv}(S_\kappa)$ of its representation in the convolution quotients field as follows:
\begin{equation}
\label{y-i-1}
P_{iv}(S_\kappa) =  \sum_{j=1}^{n} a_j S_\kappa^j \cdot \sum_{i=0}^{j-1}y_{i,0}\{1 \}^{i+1} = \sum_{i=0}^{n-1} y_{i,0}k^{i+1} P_i(S_\kappa),
\end{equation}
where $P_i(S_\kappa),\ i=0,\dots,n-1$ are polynomials in $S_\kappa$ of degree $n-i-1$:
\begin{equation}
\label{pol}
P_i(S_\kappa) = \sum_{j=i+1}^n a_j S_\kappa^{j-i-1}.
\end{equation}
For the sake of simplicity of the formulas, in what follows we suppose that the zeros $\lambda_j,\ j=1,\dots,n$ of the polynomial $P_n(S_\kappa)$ are all simple (in the case of the multiple zeros, the derivation procedure is the same, but the resulting formulas are a bit complicated).  
Then the representations of  the 
rational functions $P_i(S_\kappa)/P_n(S_\kappa),\ i=0,\dots,n-1$ in the variable $S_\kappa$  as sums of the partial fractions on the field $\mathcal{F}_{-1}$ are as follows: 
\begin{equation}
\label{partial}
\frac{P_i(S_\kappa)}{P_n(S_\kappa)} = \sum_{j=1}^n d_{ij}\frac{I}{S_\kappa - \lambda_j},\ i=0,\dots,n-1.
\end{equation}

Now we put the representations \eqref{partial} and \eqref{y-i-1} into the formula 
\eqref{eq-2-4} and employ the operational relations \eqref{op-rel} and \eqref{op-rel-1} to get the following result:
\begin{equation}
\label{y-i-2}
y_{iv}(t) = \sum_{i=0}^{n-1}y_{i,0} \tilde{y}_{i}(t),
\end{equation}
\begin{equation}
\label{y-i-3}
\tilde{y}_{0}(t) = \sum_{j=1}^n d_{0j}\, L_{\kappa,\lambda_j}(t),\ 
\tilde{y}_{i}(t) =  \left( k^{i}\, *\, \left(\sum_{j=1}^n d_{ij}\, L_{\kappa,\lambda_j}\right)\right)(t),\ i=1,\dots,n-1,
\end{equation}
where the convolution series $L_{\kappa,\lambda}$ is defined by the formula \eqref{L}. 

Summarizing the results presented above, the unique solution $y$ to the initial-value problem \eqref{eq-2-1} can be represented in the form
\begin{equation}
\label{sol2}
y(t) = y_f(t) + \sum_{i=0}^{n-1}y_{i,0} \tilde{y}_{i}(t),
\end{equation}
where the functions $y_f$ and $\tilde{y}_{i}$ are given by the formulas \eqref{y-f-4} and \eqref{y-i-3}, respectively, in terms of the convolution series  $l_{\kappa,\lambda}$ and  $L_{\kappa,\lambda}$. Moreover, our derivations implicate that 
the function $y_f$ 
is the unique solution to the problem \eqref{eq-2-1} with zero initial conditions,
and the functions $\tilde{y}_{i},\ i=0,\dots,n-1$ are the solutions to the problem \eqref{eq-2-1} for the homogeneous fractional differential equation ($f(t)\equiv 0,\ t>0$) that
fulfill the initial conditions $\tilde{y}_{i}^{(k)}(0)=\delta_{ik},\ i,k=0,\dots,n-1$. 

The inclusion $y_f\in C_{-1}^n$ follows from the formula \eqref{y-f-4}, the representation \eqref{y-f-3}, and Theorem \ref{t2}. The inclusions $\tilde{y}_{i} \in C_{-1}^n,\ i=0,\dots,n-1$ can be easily directly verified based on the representation \eqref{y-i-3} and the properties of the convolution series $L_{\kappa,\lambda}$. 

Finally, we mention that the convolution series $l_{\kappa,\lambda}$ and $L_{\kappa,\lambda}$ are the far reaching generalizations of the power series for the exponential  function ($\kappa(t) \equiv 1$, the Mikusi\'nski's operational calculus for the first order derivative) and the Mittag-Leffler function ($\kappa(t) = h_\alpha,\ 0<\alpha<1$, the Mikusi\'nski type operational calculus for the Caputo fractional derivative). In this paper, we derived some basic properties of these convolution series associated to the Mikusi\'nski type operational calculi for the GFDs with the general Sonine kernels $k \in \mathcal{S}_{-1}$. More complete theory of the convolution series will be presented elsewhere.  As we have seen, for some particular kernels $\kappa$, the convolution series $l_{\kappa,\lambda}$ and $L_{\kappa,\lambda}$ and thus the solutions  to the initial-value problem \eqref{eq-2-1} for the multi-term fractional differential equations with several GFDs can be expressed in terms of the known special functions. As an example, we refer to \cite{LucGor99} for a detailed discussion of the initial-value problems for the equations of type 
\eqref{eq-2-1} and for more general fractional differential equations with the Caputo fractional derivatives ($\kappa(t) = h_{\alpha}(t),\ 0<\alpha <1$).  

In this paper, we restricted ourselves to the mathematical properties of the GFIs and GFDs with the Sonine kernels from $\mathcal{S}_{-1}$ and analytical treatment of the initial-value problems for the fractional differential equations with these derivatives.  However, it is worth mentioning that these operators and the corresponding fractional differential equations are also potentially useful for applications.  The models with the GFDs  with the Sonine kernels from $\mathcal{S}_{-1}$ provide an essentially larger scope for fitting  the data available from the measurements compared to the models with the conventional fractional derivatives. In particular, the multi-term fractional differential equations with several GFDs are new mathematical objects that are worth for investigation both from the viewpoint of numerical methods for their numerical solution and their applications for modeling of the multi-scale time-depending processes.

\end{document}